\documentclass[12pt, reqno]{amsart}
\usepackage{pifont}
\usepackage{mathrsfs}
\usepackage{geometry}
\usepackage{titletoc}
\usepackage{stix2}
\usepackage{amscd}

\usepackage{amsmath}
\usepackage{amssymb} 
\usepackage{enumitem} 
\usepackage{mathtools} 
\usepackage[table]{xcolor} 
\usepackage[all]{xy} 
\usepackage{tikz} 
\usepackage{tikz-cd}
\usepackage{indentfirst} 
\usepackage{babel} 
\usepackage{setspace} 

\usepackage[colorlinks,linkcolor=red,anchorcolor=green,citecolor=blue]{hyperref} 
\hypersetup{linktocpage = true} 

\usepackage{rotating} 

\usepackage{ytableau} 
\usepackage{longtable} 
\newcolumntype{M}[1]{>{\centering\arraybackslash}m{#1}} 

\usepackage{fancyhdr}

\newcommand\cF{{\mathcal F}}

\newcommand\sE{{\mathscr E}}

\newcommand\sQ{{\mathscr Q}}

\newcommand{\id}{{\rm id}}
\newcommand{\rk}{{\rm rk}}

\newcommand\bp{{\bar\partial}}

\usepackage{amsthm}
\theoremstyle{plain}

\newtheorem{theorem}{Theorem}[section]
\newtheorem{thm}{Theorem}[section]
\newtheorem{lemma}[thm]{Lemma}
\newtheorem{prop}[thm]{Proposition}
\newtheorem{cor}[thm]{Corollary}
\newtheorem{defn}[thm]{Definition}

\theoremstyle{definition}
\newtheorem{example}[thm]{Example}
\newtheorem{remark}[thm]{Remark}

\newcommand{\btheorem}{\begin{thm}}
    \newcommand{\etheorem}{\end{thm}}
\newcommand{\bproposition}{\begin{prop}}
    \newcommand{\eproposition}{\end{prop}}
\newcommand{\bdefinition}{\begin{defn}}
    \newcommand{\edefinition}{\end{defn}}
\newcommand{\bcorollary}{\begin{cor}}
    \newcommand{\ecorollary}{\end{cor}}
\newcommand{\bproof}{\begin{proof}}
    \newcommand{\eproof}{\end{proof}}
\newcommand{\bremark}{\begin{remark}}
    \newcommand{\eremark}{\end{remark}}
\newcommand{\eexample}{\end{example}}
\newcommand{\bexample}{\begin{example}}

\newcommand{\elemma}{\end{lemma}}
\newcommand{\blemma}{\begin{lemma}}

\newcommand{\sq}{\sqrt{-1}}

\newcommand{\p}{\partial}

\renewcommand{\bar}{\overline}
\newcommand{\eps}{\varepsilon}

\renewcommand{\phi}{\varphi}

\newcommand{\beq}{\begin{equation}}
\newcommand{\eeq}{\end{equation}}
\newcommand{\ee}{\end{eqnarray*}}
\newcommand{\be}{\begin{eqnarray*}}

\newcommand{\bd}{\begin{enumerate}}
    \newcommand{\ed}{\end{enumerate}}

\renewcommand{\tilde}{\widetilde}

\newcommand{\qtq}[1]{\quad\mbox{#1}\quad}
\renewcommand{\bp}{\bar{\partial}}
\newcommand{\Om}{\Omega}

\newcommand{\ts}{\otimes}

\renewcommand{\>}{\rightarrow}

\newcommand{\C}{{\mathbb C}}

\newcommand{\R}{{\mathbb R}}

\renewcommand{\>}{\rightarrow}

\renewcommand{\p}{{\partial}}
\renewcommand{\bp}{{\bar{\partial}}}

\newcommand{\om}{\omega}

\renewcommand{\bar}{\overline}
\renewcommand{\tilde}{\widetilde}

\newcommand{\smo}{\sqrt{-1}}

\renewcommand{\id}{\mathrm{Id}}

\newcommand{\Herm}{\mathrm{Herm}}
\newcommand{\tr}{\mathrm{tr}}
\newcommand{\nm}[1]{\left\Vert #1\right\Vert}
\setlist[itemize]{leftmargin=*}
\setlist[enumerate]{leftmargin=*}

\numberwithin{equation}{section} 

\makeatletter

\usepackage{fancyhdr}
\renewcommand{\leftmark}{{\scriptsize Existence of twisted Hermitian-Einstein metrics on unstable vector bundles}}

\title{Existence of twisted Hermitian-Einstein metrics on unstable vector bundles}

\author{Mingwei Wang}
\author{Xiaokui Yang}
\author{Shing-Tung Yau}

\address{Mingwei Wang, Qiuzhen College, Tsinghua University, Beijing, 100084, China}
\email{wangmw21@mails.tsinghua.edu.cn}

\address{Xiaokui Yang, Department of Mathematics and Yau Mathematical Sciences Center, Tsinghua University, Beijing, 100084, China}
\email{xkyang@mail.tsinghua.edu.cn}

\address{Shing-Tung Yau, Yau Mathematical Sciences Center and  Qiuzhen College, Tsinghua University, Beijing, 100084, China}
\email{styau@mail.tsinghua.edu.cn}

\begin{document}

    \begin{abstract} In this paper, we demonstrate that twisted Hermitian-Einstein metrics on holomorphic vector bundles exist without obstruction. More precisely, for an arbitrary  holomorphic vector bundle $E$ over a compact K\"ahler manifold $(M,\omega_g)$, we prove that the  twisted Hermitian-Einstein equation
    $$
    	\Lambda_{\omega_g}\left(\sqrt{-1}R^h\right) = \lambda h + P
    $$
    admits a unique smooth  solution $h$,  provided that  $P\in\Gamma(M,E^*\otimes\bar{E}^*)$ is  positive-definite and  $\lambda<\lambda_E^-$. The constant  $\lambda_E^-$ is intrinsically associated with the stability constant  of $E$.
This result extends the classical Donaldson–Uhlenbeck–Yau (DUY) theorem for stable  bundles and, in the limit $P\>0$, gives a new proof of the DUY theorem. As an application, we obtain an intrinsic Chern number inequality for \emph{unstable vector bundles}: $$\int_M \left((r-1)c_1(E)^2 - 2rc_2(E)\right) \wedge \omega_g^{n-2} \leq \Bigl\lfloor \frac{r^2}{4} \Bigr\rfloor \frac{(\lambda_E^+-\lambda_E^-)^2}{4\pi^2 n^2} \int_M \omega_g^n.$$
    \end{abstract}

    \maketitle {
        \setcounter{tocdepth}{1}

    {\small{    \begin{spacing}{1.1} \tableofcontents %
                \dottedcontents{section}[1.8cm]{}{3em}{5pt} %
\end{spacing} }} }

\vspace*{-1\baselineskip} 
    \section{Introduction}

This paper continues the investigation initiated in \cite{WYY26+} and \cite{FWYY26+},  where  the prescribed Hermitian-Yang-Mills tensor equation for holomorphic (Higgs) vector bundles on compact Hermitian manifolds was introduced:
\beq \Lambda_{\om_g}\left(\sqrt{-1}\, R^h\right) =P, \label{PHYM}\eeq
which can be viewed as a vector bundle analogue of the Calabi–Yau theorem (\cite{Cal57,Yau78}). It is well-known that the Calabi-Yau theorem asserts that: on a compact K\"ahler manifold $(M,\omega_g)$,  for any closed real $(1,1)$ form $\Om$ representing $2\pi c_1(M)$,  there exists a unique K\"ahler metric $\omega\in [\omega_g]$ such that its  Ricci curvature satsifies \beq \mathrm{Ric}(\omega)=\Om.\label{CY}\eeq

Let's recall the classical Donaldson--Uhlenbeck--Yau theorem (\cite{Don85, UY86, Don87}) for stable vector bundles (see also the Hermitian analogue in \cite{LY86}):  \setcounter{theorem}{0}
\renewcommand{\thetheorem}{\Alph{theorem}}\begin{theorem} \label{DUY} For  a stable holomorphic vector bundle $E$ over a compact K\"ahler manifold $(M,\omega_g)$,  there exists a unique Hermitian-Einstein metric $h$ on $E$ up to scaling, satisfying
    \beq \Lambda_{\omega_g} \left(\sq R^h\right)=\lambda_E\cdot h, \label{HE} \eeq
    where $R^h\in\Gamma(M,\Lambda^{1,1}T^*M\ts E^*\ts \bar E^*)$ is the Chern curvature tensor of $(E,h)$. 
 \end{theorem}
\noindent Slope stability is a fundamental algebraic condition:  for every non-zero proper coherent subsheaf $\mathcal{F} \subset E$, the following inequality holds:
\beq 
\frac{\deg_{\omega_g}(\mathcal{F})}{\operatorname{rk}(\mathcal{F})} < \frac{\deg_{\omega_g}(E)}{\operatorname{rk}(E)}.
\label{stable condition}
\eeq 
Here the degree is  defined by
\beq 
\deg_{\omega_g}(\mathcal{F}) = 2\pi n \int_M c^{\mathrm{BC}}_1(\mathcal{F}) \wedge \omega_g^{n-1},
\eeq 
where $c^{\mathrm{BC}}_1(\cF)$ denotes the first Bott--Chern class in $H^{1,1}_{\mathrm{BC}}(M, \mathbb{C})$, coinciding with the usual first Chern class $c_1(\cF)$ when $M$ is K\"ahler. We refer to \cite{Kob87} for more details.  
\noindent For a coherent subsheaf  or quotient sheaf  $\sQ$ of $E$, the \emph{stability constant} $\lambda_\sQ$ is defined as  \beq
\lambda_\sQ:=\frac{2\pi n\int_M c^{\mathrm{BC}}_1(\sQ)\wedge
	\omega_g^{n-1}}{\mathrm{\rk}(\sQ)\int_M\omega_g^n}=\frac{\deg_{\omega_g}(\sQ)}{\rk(\sQ)
	\int_M\omega_g^n}.\eeq 
The following extremal constants are well-known:
 \beq \lambda_E^+ := \sup\left\{ \left. \lambda_\mathscr{S}\ \right|\  \mathscr{S} \text{ is a coherent subsheaf of } E, \mathrm{rk}(\mathscr{S}) > 0 \right\}, \label{lambda+}\eeq
and
  \beq \lambda_E^- := \inf\left\{ \left. \lambda_\sQ\ \right|\  \mathscr{Q} \text{ is a coherent quotient sheaf of } E, \mathrm{rk}(\mathscr{Q}) > 0 \right\}. \label{lambda-}\eeq
One immediately obtains the inequality
\beq \lambda_E^-\leq \lambda_E\leq \lambda_E^+.\eeq 
Moreover, $E$ is semi-stable if and only if $\lambda_E^-=\lambda_E$ or $\lambda_E^+=\lambda_E$; in this case,  \beq  \lambda_E^-= \lambda_E=\lambda_E^+.\eeq

\vskip 1\baselineskip

The Calabi-Yau equation \eqref{CY} and the Hermitian-Einstein equation \eqref{HE} have been studied extensively over the past half-century, yielding numerous important applications and emerging as central topics in  differential geometry, algebraic geometry, mathematical physics and etc.. For a comprehensive account of this topic, we refer to \cite{NS65, Siu87, Sim88, Sim92, Kol98, DW04, FY08, TW10, Jac14, LZ15, CS16, STW17, Sze18, CSW18, Fu19, CS20, DS21, Li21, MPS25} and references therein.\\

The prescribed Hermitian-Yang-Mills tensor equation \eqref{PHYM} interpolates between the Calabi-Yau equation \eqref{CY} and the Hermitian-Einstein equation \eqref{HE}. By inheriting the geometric and analytic features of both, it serves as a pivotal bridge connecting the geometry of vector bundles to the broader realm of K\"ahler and Hermitian geometry.\\

    The main result of this paper establishes  the existence of twisted Hermitian-Einstein metrics on  arbitrary  holomorphic vector bundles.

    \btheorem
    \label{main}
    Let $ (M,\omega_g) $ be a compact K\"ahler (or Gauduchon) manifold and $ E $ be an arbitrary holomorphic vector bundle over $ M $.  If we define 
     \beq \mu_E^-: = \sup
    \left\{\ \lambda\ \left|\  \text{ there exists a Hermitian metric }
    h \text{ such that } \Lambda_{\omega_g}\left(\sq R^h\right) >
    \lambda h \right.\right\}, \label{mu} \eeq
     then \beq \lambda_E^-=\mu_E^-.\eeq  Moreover,   for every $\lambda<\lambda_E^-$ and every positive-definite Hermitian tensor $P\in \Gamma(M,E^*\ts \bar E^* )$, there exists a unique Hermitian metric $h$ on $E$ satisfying the twisted Hermitian-Einstein equation
    \beq \Lambda_{\omega_g}\left(\sq R^h\right)=\lambda\cdot h+ P.\label{tHE}\eeq
    \etheorem
    \noindent 

\noindent
    One of the key ingredients in the proof of Theorem \ref{main} is the following existence result for the twisted Hermitian-Einstein equation:

\btheorem
\label{main2}
Let $ (M,\omega_g) $ be a compact Hermitian manifold and $ E $ be a holomorphic vector bundle over $ M $. For given $ \lambda \in \mathbb{R} $, if there exists a smooth Hermitian metric $ h_0 $ on $ E $ such that
\beq \Lambda_{\omega_g}\left(\sq R^{h_0}\right) - \lambda h_0  > 0, \label{curvaturecondition} \eeq
then for any positive-definite Hermitian tensor $ P \in \Gamma(M,E^*\ts \bar E^*) $, there exists a unique smooth Hermitian metirc $ h $ on $ E $ satisfying the twisted Hermitian-Einstein equation:
\beq \Lambda_{\omega_g}\left(\sq R^h\right)=\lambda h+ P.\label{THE} \eeq
\etheorem

\noindent The proof of Theorem~\ref{main2} proceeds as follows. The uniform a priori  $C^0$-estimate and the uniqueness for the twisted Hermitian-Einstein equation \eqref{THE} follow from the comparison principle below.
\btheorem\label{main5} Let $ (M,\omega_g)
$ be a compact Hermitian manifold and $ E $ be a holomorphic vector
bundle over $ M $. If there are  Hermitian metrics $ h$ and  $h_0 $
on $E$ and a constant $\lambda\in \R$ satisfying
$\Lambda_{\omega_g}\left(\sq R^{h_0}\right)-\lambda h_0>0$, and
\beq  \Lambda_{\omega_g}\left(\sq R^{h}\right) -
\lambda h \leq \Lambda_{\omega_g}\left(\sq R^{h_0}\right) - \lambda
h_0, \label{Hermitiantensor} \eeq as Hermitian tensors in $ \Gamma(M,E^*\otimes \bar E{}^*)
$. Then $ h \leq h_0 $. \etheorem
\noindent 
Note that inequality \eqref{Hermitiantensor} holds as Hermitian tensors in $\Gamma(M,E^*\otimes  \bar E^*)$,  \textbf{but not  as scaling invariant endomorphisms} in  $\Gamma(M,E^*\otimes E)$.\\

The existence of solutions to \eqref{THE} is established by  an implicit function theorem argument
and the following weak compactness result for twisted Hermitian-Einstein equations.
\btheorem\label{main4} Let $ (M,\omega_g) $ be a compact Hermitian manifold and $ (E, h_0)$ be a Hermitian vector bundle over $M$. Fix $\lambda\in \R$.  Suppose that the following conditions are satisfied  \bd \item
$\{ P_m \} \subset \mathrm{Herm}(E)$
converges smoothly to some $P \in \mathrm{Herm}(E)$;
\item $\{\lambda_m\}\subset \R$ converges to $\lambda$;
\item  $\{h_m\} \subset \mathrm{Herm}^+(E)$ satisfies
\beq 
\Lambda_{\omega_g} \left(\sqrt{-1} R^{h_m}\right) = \lambda_m h_m + P_m;\eeq
\item  there exists a uniform constant $C>0$ such that
\begin{equation}
C^{-1} h_0 \leq h_m \leq C h_0,
\end{equation}
\ed  holds for all $m$.
Then, there exists a subsequence $\{h_{m_i}\}$ which converges smoothly to a smooth Hermitian metric $h$ on $E$, and this limit satisfies
\begin{equation}
\Lambda_{\omega_g}\left(\sq  R^{h} \right)= \lambda h + P.
\end{equation}
\etheorem

\noindent We briefly outline the proof of Theorem \ref{main}. The inequality $\mu_E^- \leq \lambda_E^-$ follows easily
from the definitions. We proceed by contradiction and assume $\mu_E^- <
\lambda_E^-$. By Theorem \ref{main2}, for a sequence $\mu_m \nearrow
\mu_E^-$, there exist Hermitian metrics $\{h_m\}$ solving \beq
\Lambda_{\omega_g}\left(\sq R^{h_m}\right) = \mu_m h_m+P. \eeq The comparison principle Theorem \ref{main5}
implies that $\{h_m\}$ are uniformly bounded from below. \emph{Case
1:} If $\{h_m\}$ is unbounded, a technique adapted from \cite{UY86}
allows us to construct a coherent quotient sheaf $\sQ$ with
$\lambda_{\sQ}$ satisfying $\lambda_E^- \leq \lambda_{\sQ} \leq \mu_E^-$,
contradicting $\mu_E^- < \lambda_E^-$. \emph{Case 2:} If $\{h_m\}$ is
uniformly bounded, Theorem
\ref{main4} yields a limit metric $h$ solving
\begin{equation}\label{eq:contradiction_limit}
\Lambda_{\omega_g}\bigl(\sqrt{-1}\, R^h\bigr) - (\mu_E^-+\varepsilon)h = P-\varepsilon h>0,
\end{equation}
for some $\varepsilon > 0$. This contradicts the maximality of $\mu_E^-$.
Thus $\mu_E^- = \lambda_E^-$, and Theorem \ref{main} follows immediately from Theorem \ref{main2}.\\

As an application of Theorem \ref{main} and Theorem \ref{main2}, we obtain a precise quantitative approximation for the Hermitian-Yang-Mills tensors of arbitrary vector bundles: 
\btheorem\label{main7}
Let $ (M,\omega_g) $ be a compact Gauduchon manifold and $ E $ be a holomorphic vector bundle over $ M $. For any $ \eps > 0 $, there exists a Hermitian metric $ h_\eps $ on $ E $ such that 
\beq \left(\lambda_E^-- \eps\right) \cdot h_\eps \leq \Lambda_{\omega_g}\left(\smo R^{h_\eps}\right) \leq \left(\lambda_E^++ \eps\right) \cdot h_\eps. \eeq
\etheorem

As a consequence of Theorem~\ref{main7}, we establish an intrinsic Chern number inequality for  arbitrary holomorphic vector bundles, generalising  the classical inequalities known for semi-stable vector bundles.
\btheorem \label{main6}Let $ (M,\omega_g) $ be a compact K\"ahler manifold and $ E $ be an arbitrary holomorphic vector bundle over $ M $ of rank $r$. The following Chern number inequality holds:
\beq\int_M \left((r-1)c_1(E)^2 - 2rc_2(E)\right) \wedge \omega_g^{n-2} \leq \Bigl\lfloor \frac{r^2}{4} \Bigr\rfloor \frac{(\lambda_E^+-\lambda_E^-)^2}{4\pi^2 n^2} \int_M \omega_g^n. \eeq
\etheorem

Moreover, we obtain a new proof of the classical Donaldson-Uhlenbeck-Yau theorem on the existence of  Hermitian-Einstein metrics on stable vector bundles (\cite{Don85,UY86,Don87}, \cite{LY86}).

\btheorem\label{main3}
\label{alternativesemistable}\label{alternativeUBY}
Let $ (M,\omega_g) $ be a compact K\"ahler (or Gauduchon) manifold and $ E $ be a  holomorphic vector bundle over $ M $.
 If $E$ is semi-stable, then for any $ P \in \mathrm{Herm}^+(E) $ and for any $ \eps > 0 $, there exists a unique smooth Hermitian metric $ h_\eps $ on $E$ such that
\beq \Lambda_{\omega_g}\left(\smo R^{h_\eps}\right) = (\lambda_E - \eps)\cdot h_\eps + \eps P,  \label{HEt}\eeq
and $ E $ is approximately Hermitian-Einstein.
Moreover, if $E$ is stable, then a subsequence $\{h_{\eps_i}\}$ converges smoothly to a Hermitian-Einstein metric $h$ on $E$ as $\eps_i\>0$.
\etheorem
\noindent The argument for Theorem~\ref{main3} is not a straightforward adaptation
of Theorem~\ref{main} or the methods in \cite{Don85, UY86, Don87, LY86}. Indeed, the family of equations \eqref{HEt} is established by Theorem~\ref{main} and the uniform $C^0$-estimate
	\begin{equation}
	C \cdot h_0 \le h_\varepsilon
	\end{equation}
	is obtained by detailed calculations together with the Moser iteration,
	rather than via the comparison principle established in this paper.

\bremark The main results of the paper  also hold for Higgs bundles over compact Hermitian manifolds,
by adapting the frames constructed in \cite{FWYY26+}.
Further applications are obtained in \cite{XYY26+} via the associated parabolic flows of the form
\begin{equation}
\label{GHYMF}
\frac{\partial h}{\partial t} = -\Lambda_{\omega_g}\left(\sqrt{-1}\, R^h\right) + \mathcal{F}_0(h),
\end{equation}
where $\mathcal{F}_0$ is a general map taking values in $\mathrm{Herm}(E)$ and $R^{h}$ represents the curvature of an affine connection; in particular, it includes the curvature associated with a Higgs bundle. Specific instances include:
\begin{enumerate}[label=(\roman*)]
	\item $\mathcal{F}_0(h) = P$ for some fixed $P \in \mathrm{Herm}(E)$;
	\item $\mathcal{F}_0(h) = \lambda h+P$;
	\item $\mathcal{F}_0(h) = \lambda h + P + (\mathrm{tr}_{h_0} h) h_0$ for some fixed Hermitian metric $h_0$ on $E$.
\end{enumerate}
Moreover, the positivity assumption on $P$ can be weakened using flow techniques.
\eremark

\noindent\textbf{Acknowledgements}. The second named author would like to thank Bing-Long Chen, Jixiang Fu, Kefeng Liu, Song Sun and Valentino Tosatti for inspiring  discussions.

\vskip 2\baselineskip

\section{A comparison principle and proof of Theorem \ref{main5}}

We recall some background material  for the reader's convenience.
Let $(M,\omega_g)$ be a compact Hermitian manifold of dimension $n$.
We say that $(M,\omega_g)$ is a \emph{Gauduchon manifold} if
$
\partial\bar{\partial}\,\omega_g^{\,n-1}=0.
$
It is well known that, for any Hermitian metric $\omega$ on $M$,
there exists a smooth function $f$ such that
$
\omega_g = e^{f}\omega
$
satisfies $\partial\bar{\partial}\,\omega_g^{\,n-1}=0$ (\cite{Gau84}). Let $E$ be a holomorphic vector bundle over $M$ of rank $r$. We use the following notations:

\begin{itemize}
    \item $\mathrm{Herm}(E)$: the space of Hermitian tensors in
    $\Gamma(M, E^*\otimes \bar{E}^*)$;

    \item $\Herm^+(E)$: the subspace consisting of
    positive-definite Hermitian tensors;

    \item $\Herm^{\ge 0}(E)$: the subspace consisting of
    non-negative Hermitian tensors.
\end{itemize}
Given a smooth Hermitian metric $h$ on $E$, the Chern connection of
$(E,h)$ is denoted by $\nabla^{h}$. We shall use the natural
decomposition $ \nabla^{h}=\p^{h}+\bp$, where $\p^{h}$ is the
$(1,0)$ part and $\bp$ is the $(0,1)$ part. Let $R^{h}$ be the Chern
curvature tensor of $(E,h)$. In local holomorphic coordinates
$\{z^i\}$ of $M$ and local holomorphic basis $\{e_\alpha\}$ of $E$,
\beq R^h=R_{i\bar j\alpha\bar\beta}dz^i\wedge d\bar z^j\ts
e^\alpha\ts\bar{e}^\beta \in\Gamma\left(M,\wedge^{1,1}T^*M\ts
E^*\ts\bar{E}^*\right). \eeq  The \textbf{Hermitian-Yang-Mills tensor} $ S^{h}\in
\Gamma(M,E^*\ts \bar E^*)$ of $(E,h)$  is defined as \beq S^{h}
:=\Lambda_{\omega_g}\left(\smo R^{h}\right)=\left(g^{i\bar
    j}R^{h}_{i\bar j\alpha\bar\beta}\right) e^\alpha\ts \bar e^\beta \in
\Gamma(M,E^*\ts \bar E^*).
\eeq  By using the Hermitian metric $h$, there
are natural lifts of $R^h$ and $S^h$: \beq\Theta^h=R^h\cdot
h^{-1}=R_{i\bar j\alpha}^{\beta}dz^i\wedge d\bar z^j\ts
e^\alpha\ts{e}_\beta\in \Gamma(M,\Lambda^{1,1}T^*M\ts E^*\ts E),\eeq
and \beq   K^h=S^h\cdot h^{-1}=\left(g^{i\bar j}R_{i\bar
    j\alpha}^{\beta}\right) e^\alpha\ts  e_\beta\in \Gamma(M,E^*\ts E).
\eeq

\noindent The following transformation formula for the
Hermitian-Yang-Mills curvature tensors is essentially well-known
(e.g. \cite[Proposition~3.6]{WYY26+})
\blemma\label{lem Theta formula} Let $(M,\omega_g)$ be
a compact  Hermitian manifold and $E\>M$ be a holomorphic vector
bundle. If $h$ and $h_0$ are smooth Hermitian metrics on $E$, then
\beq \Lambda_{\omega_g}\left(\sq  \Theta^{h}\right)-
\Lambda_{\omega_g}\left(\sq  \Theta^{h_0}\right)=
\Lambda_{\omega_g}\smo\bp\left( (\p^{h_0} H) \cdot H^{-1}
\right)\label{conformalchange1} \eeq as tensors in $\Gamma(M,E^*\ts
E)$ where $H=h\cdot h_0^{-1}$.
\elemma

\noindent The following linear algebraic result is  established in \cite[Lemma~2.2]{WYY26+}:
\blemma \label{linearalgebra}\bd
\item $H=h\cdot h_0^{-1}\in\Gamma(M,E^*\ts E)$    is  $h_0$-Hermitian.

\item If $P\in \Gamma(M,E^*\ts E)$ is $h_0$-Hermitian,  for any $ A, B \in \Gamma(M,E^*\otimes E) $,
\beq h_0(P\cdot A, B)=h_0(A, P\cdot B), \quad  h_0(A \cdot P, B) = h_0(A, B\cdot P). \eeq

\item If $ A, B \in \Gamma(M,E^*\otimes E) $ are $h_0$-Hermitian and $ A\geq 0, B \geq 0 $, then
\beq h_0(A,B) \geq 0. \eeq Moreover,  for  any section $ C \in \Gamma(M,E^* \otimes E) $, one has
\beq h_0(A \cdot C \cdot B, C ) \geq 0. \eeq

\ed
\elemma

\noindent The main result of this section is the following
comparison principle: \btheorem\label{comparison} Let $ (M,\omega_g)
$ be a compact Hermitian manifold and $ E $ be a holomorphic vector
bundle over $ M $. If there are  Hermitian metrics $ h$ and  $h_0 $
on $E$ and a constant $\lambda\in \R$ satisfying
$\Lambda_{\omega_g}\left(\sq R^{h_0}\right)-\lambda h_0>0$, and
\beq \label{Conditioncom} \Lambda_{\omega_g}\left(\sq R^{h}\right) -
\lambda h \leq \Lambda_{\omega_g}\left(\sq R^{h_0}\right) - \lambda
h_0, \eeq as Hermitian tensors in $ \Gamma(M,E^*\otimes \bar E{}^*)
$. Then $ h \leq h_0 $. \etheorem

\bproof
Let $ H = h \cdot h_0^{-1} \in \Gamma(M,E^*\otimes E) $.  Define
\beq \label{DefofOm}\Omega^\lambda:= \Lambda_{\omega_g}\left(\sq \Theta^{h_0}\right) - \lambda\mathrm{Id}_E \in \Gamma(M,E^*\ts  E), \eeq
and 
\beq \Phi^\lambda:= \left(\Lambda_{\omega_g}\left(\sq \Theta^{h}\right) -  \lambda\mathrm{Id}_E\right) \cdot H \in \Gamma(M,E^*\ts  E). \eeq 
In particular, the conditions in Theorem \ref{comparison} imply that $ \Omega^{\lambda} > 0 $ and \beq \Omega^{\lambda} \leq \Phi^{\lambda}\eeq  with respect to $h_0$. On the other hand, by \eqref{conformalchange1}, we have
\beq \Phi^{\lambda}\cdot H^{-1}-\Om^\lambda= \Lambda_{\omega_g} \left(\smo \Theta^h\right) -\Lambda_{\omega_g} \left(\smo \Theta^{h_0}\right) = \Lambda_{\omega_g}\smo \bp\left(\p^{h_0}H \cdot H^{-1}\right).\label{curvaturedifference} \eeq
This is equivalent to the following form:
\beq \Omega^{\lambda} \cdot (H^{-1} - \mathrm{Id}_E) = \Lambda_{\omega_g}\smo \bp\left(\p^{h_0}H \cdot H^{-1}\right) + (\Omega^{\lambda} - \Phi^{\lambda}) \cdot H^{-1}.\eeq 
The  proof of $H\leq \mathrm{Id}_E$ then follows by a strategy similar to that used in the proof of \cite[Theorem~1.3]{FWYY26+}. \eproof

As an application of Theorem \ref{comparison}, we obtain a version of  uniform $C^0$-estimates:

\bproposition Let $E$ be a holomorphic vector bundle over  a compact Hermitian manifold $(M,\omega_g)$. Fix a constant $\lambda\in \R$. Suppose that $h_0$ and $h$ are two Hermitian metrics on $E$, and \beq \Lambda_{\omega_g}\left(\sq R^{{h_0}}\right)-\lambda h_0>0. \label{initialpositivity}\eeq
If there exists a constant  $C_1\geq 1 $ such that
\beq C_1^{-1} h_0\leq  \Lambda_{\omega_g}\left(\smo R^{{h}}\right)-\lambda h\leq C_1 h_0 \label{curvatureuniformbound0} \eeq
as Hermitian tensors in $\Gamma(M,E^*\ts  \bar E^*)$, then
there exists $C_2=C_2(\omega_g,  h_0, C_1)$ such that
\beq C_2^{-1}h_0 \leq h \leq C_2 h_0. \eeq
\eproposition
\bproof
Since $M$ is compact, by \eqref{initialpositivity} and \eqref{curvatureuniformbound0}, there exists $C_2=C_2(\omega_g,  h_0, C_1)$ such that
\beq C_2^{-1}\left(\Lambda_{\omega_g}\left(\sq R^{h_0}\right) - \lambda h_0\right)\leq \Lambda_{\omega_g}\left(\sq R^{h}\right) - \lambda h\leq C_2\left(\Lambda_{\omega_g}\left(\sq R^{h_0}\right) - \lambda h_0\right).\eeq
Applying Theorem~\ref{comparison} with a suitable scaling yields
\beq C_2^{-1}h_0 \leq h \leq C_2 h_0. \eeq
This completes the proof.
\eproof

\vskip 2\baselineskip

\section{A priori estimates}

In this section, we establish the following a priori estimates:
\btheorem\label{C1estimate} Let $E$ be a holomorphic vector bundle over  a compact Hermitian manifold $(M,\omega_g)$. Fix a constant $\lambda\in \R$.  Suppose that $h_0$ and $h$ are two Hermitian metrics on $E$.  We set $H=h\cdot h_0^{-1}$ and
\beq \Phi^\lambda:=  \left(\Lambda_{\omega_g}\left(\smo R^{{h}}\right)\right)\cdot h_0^{-1}-\lambda H \in \Gamma(M,E^*\ts E).\eeq
If there exist constants  $C$, $c_1$ and $c_2$ such that $|\lambda|<c_1$, $|\Phi^\lambda|_{C^1(M,\omega_g, h_0)}<c_2$
and \beq C^{-1} h_0 \leq  h\leq C h_0, \label{curvatureuniformbound} \eeq
then  the following  $C^1$-estimate holds:
\beq  |H|_{C^1(M,\omega_g, h_0)}\leq C_1,\eeq
where the constant $C_1$ depends on $M,\omega_g,  h_0, C, c_1$ and $c_2$.
\etheorem

Let $ \sE: = T^{*1,0}M \otimes E^* \otimes E $ be the Hermitian holomorphic vector bundle with the Hermitian metric induced by $ g $ on $ T^{*1,0}M $ and $ h $ on $ E $.  Let  $ T $ be the tensor
\beq T := \p^{h_0}H \cdot H^{-1}  \in \Omega^{1,0}(M,E^*\otimes E) \simeq \Gamma(M,\sE). \eeq

\bproposition
\label{DeltaS} Let $E$ be a holomorphic vector bundle over  a compact Hermitian manifold $(M,\omega_g)$. Fix a constant $\lambda\in \R$.  Suppose that $h_0$ and $h$ are two Hermitian metrics on $E$.   If there exist constants  $C$, $c_1$ and $c_2$ such that $|\lambda|<c_1$, $|\Phi^\lambda|_{C^1(M,\omega_g, h_0)}<c_2$
and 
\beq C^{-1}h_0\leq h\leq Ch_0, \eeq
then
one has
\beq \Delta_{\C} |T|_h^2\geq - C_{2}\left(|T|^2_h+1\right), \eeq
where $\Delta_\C =\mathrm{tr}_{\omega_g}\sq \p\bp$ on functions, and $C_{2}$ depends on $M,\omega_g,  h_0, C, c_1$ and$ c_2$.
\eproposition

\bproof We refine the arguments from \cite{WYY26+} and \cite{FWYY26+}.  We set $$ \Phi:= \left(\Lambda_{\omega_g}\left(\smo R^{{h}}\right)\right)\cdot h_0^{-1}$$ and so
\beq \Phi^\lambda=\Phi-\lambda H.\eeq
The Bochner-Kodaira formula establishes that
\beq \Delta_{\C} |T|_h^2 = |\p_\sE T|_h^2 + |\bar\p T|_h^2 + 2\mathrm{Re}\left(\left\langle \mathrm{tr}_{\omega_g}\left(\smo \p_\sE\bar\p_\sE T\right), T \right\rangle_\sE \right) - \mathrm{Ric}^\sE(T,T),\label{c30} \eeq
where $\mathrm{Ric}^\sE$ is the Hermitian-Yang-Mills tensor of $\sE$. In \cite[Proposition~4.6]{FWYY26+}, it is established that
\beq \Delta_{\C} |T|_h^2 \geq -C_{3}\left(|\p^{h_0}\Phi|_h|T|_h+ |\Phi|_h |T|^2_h+|T|^2_h+1\right), \label{C1maininequality} \eeq
where $C_{3}$  depends on $M, \omega_g, h_0$ and $C$.
Since $\Phi=\Phi^\lambda+\lambda H$, one has
\beq
|\Phi|_h\leq |\Phi^\lambda|_h+|\lambda||H|_h, \quad
|\p^{h_0}\Phi|_h\leq |\p^{h_0}\Phi^\lambda|_h+|\lambda||\p^{h_0}H|_h.
\eeq
Since
$C^{-1}\mathrm{Id}_E\leq H\leq C\mathrm{Id}_E$ and $ T=\p^{h_0}H\cdot H^{-1}$,
one has
\beq
|\Phi|_h\leq |\Phi^\lambda|_h+C_4, \quad
|\p^{h_0}\Phi|_h\leq |\p^{h_0}\Phi^\lambda|_h+C_5|T|_h,\label{c2}
\eeq
where $C_4$ and $C_5$ depend on $M,\om_g,C,h_0$ and $c_1$. By \eqref{c2} and \eqref{C1maininequality}, one has
\beq
\Delta_{\C} |T|_h^2 \geq-C_{6}\left(|\p^{h_0}\Phi^\lambda|_h|T|_h+|\Phi^\lambda|_h|T|^2_h +|T|^2_h+1\right),
\eeq
where $C_{6}$  depends on $M, \omega_g, h_0, C$ and $c_1$. By   the  Cauchy-Schwarz inequality,  we establish   \beq \Delta_{\C} |T|_h^2\geq - C_{7}\left(|T|^2_h+1\right) \eeq
where $ C_{7}$ depends on $M,\omega_g, C, h_0,c_1$ and $c_2$.
\eproof

\noindent The following result is exactly \cite[Proposition~4.4]{FWYY26+}.

\bproposition
\label{W12estimate}   Let $E$ be a holomorphic vector bundle over  a compact Hermitian manifold $(M,\omega_g)$. Suppose that $h_0$ and $h$ are two Hermitian metrics on $E$, and $H=h\cdot h_0^{-1}$.
Suppose that there exists a constant $C$ such that
\beq C^{-1}h_0\leq h\leq Ch_0. \eeq
Then the following estimates hold:
\beq C^{-2} h_0(\p^{h_0} H, \p^{h_0} H) \leq h(\p^{h_0} H\cdot H^{-1}, \p^{h_0} H\cdot H^{-1}) \leq C^2 h_0(\p^{h_0} H, \p^{h_0} H), \label{c15}\eeq   and
\beq \mathrm{tr}_E\left(\Lambda_{\omega_g}\left(\sq \p^{h_0} H\cdot H^{-1}\cdot \bp H\right)\right)\geq C^{-1} h\left(\p^{h_0}H\cdot H^{-1}, \p^{h_0}H\cdot H^{-1}\right). \label{estimate2}\eeq
\eproposition

\bproposition\label{trEH}  Let $E$ be a holomorphic vector bundle over  a compact Hermitian manifold $(M,\omega_g)$. Fix a constant $\lambda\in \R$. Suppose that $h_0$ and $h$ are two Hermitian metrics on $E$.
Suppose that there exists a constant $C$ such that
\beq C^{-1}h_0\leq h\leq Ch_0. \eeq
Then the following estimate holds:
\beq \Delta_{\C}\mathrm{tr}_E H  \geq  C^{-1}|\p^{h_0}H \cdot H^{-1}|^2_h +\tr_E(\Om^{\lambda}\cdot H) -\mathrm{tr}_E\Phi^{\lambda},\eeq
where  $\Omega^\lambda:= \Lambda_{\omega_g}\left(\sq \Theta^{h_0}\right) - \lambda\mathrm{Id}_E \in \Gamma(M,E^*\ts  E)$.
\eproposition
\bproof  It is known that
\beq \Delta_\C \mathrm{tr}_E H=\Lambda_{\omega_g} \sq \p\bp h_0(H, \mathrm{Id}_E)=-\sq \Lambda_{\omega_g} h_0(\bp \p^{h_0} H, \mathrm{Id}_E)=-\tr_E(\Lambda_{\omega_g}\smo\bar\p\p^{h_0} H). \eeq
 By \eqref{conformalchange1},
\beq \Lambda_{\omega_g} \left(\smo \Theta^h\right) -\Lambda_{\omega_g} \left(\smo \Theta^{h_0}\right) = \Lambda_{\omega_g}\smo \bp\left(\p^{h_0}H \cdot H^{-1}\right).\label{curvaturedifference0} \eeq
Multiplying both sides by $H$ yields:
\beq \Phi^\lambda-\Om^{\lambda}\cdot H=\Lambda_{\om_g}\smo \bp(\p^{h_0}H \cdot H^{-1}) \cdot H. \eeq
Hence,
\be \Delta_{\C}\mathrm{tr}_E H & =& -\tr_E(\Lambda_{\omega_g}\smo\bar\p\p^{h_0} H) \\ &=& \mathrm{tr}_E\left(\Lambda_{\omega_g}\smo(\p^{h_0}H \cdot H^{-1} \cdot \bar\p H)\right) - \mathrm{tr}_E(\Lambda_{\omega_g}\smo\bp (\p^{h_0}H \cdot H^{-1}) \cdot H )  \\
& \geq &C^{-1}|\p^{h_0}H \cdot H^{-1}|^2_h+\tr_E(\Om^{\lambda}\cdot H) -\mathrm{tr}_E\Phi^\lambda , \ee
where the last inequality follows from \eqref{estimate2}.
This completes the proof.
\eproof

\bproof[Proof of Theorem \ref{C1estimate}] By Proposition \ref{trEH}, we have
\beq \Delta_{\C}\mathrm{tr}_E H  \geq  C^{-1}|T|^2_h+\tr_E(\Om^{\lambda}\cdot H) -\mathrm{tr}_E\Phi^{\lambda}. \eeq
Since $|\Phi^\lambda|_{C^1(M,\omega_g, h_0)}<c_2$,
\beq |\mathrm{tr}_E\Phi^{\lambda} |\leq C_2\eeq 
where $C_2=\rk(E) c_2$.  Moreover, since $$\Omega^\lambda= \Lambda_{\omega_g}\left(\sq \Theta^{h_0}\right) - \lambda\mathrm{Id}_E, \quad C^{-1}H\leq \mathrm{Id}_E\leq C\mathrm{Id}_E, $$ one has 
\beq |\tr_E(\Om^{\lambda}\cdot H)|\leq \left|\tr_E\left(\Lambda_{\omega_g}\left(\sq \Theta^{h_0}\right)\cdot H\right)\right| +|\lambda||\tr_EH|  \leq C_3, \eeq 
where $C_3$ depends on $M,\omega_g, h_0, C$ and $c_1$. Hence, 
\beq \Delta_{\C}\mathrm{tr}_E H  \geq  C^{-1}|T|^2_h-C_4 \eeq
where $C_4=C_2+C_3$. By Proposition \ref{DeltaS},
\beq \Delta_{\C} |T|_h^2 \geq -C_{5}(|T|_h^2 + 1), \eeq
where  $C_5$ depend on $M,\omega_g,  h_0, C, c_1 $ and $c_2$. For a large $ L > 0 $ such that $LC^{-1} - C_{5}>1$, one has
\beq \Delta_{\C} (|T|_h^2 + L\mathrm{tr}_E H) \geq (LC^{-1} - C_{5})|T|_h^2 - (LC_4+C_5). \eeq
At a maximal point $ p \in M $ of $ (|T|_h^2 + L\mathrm{tr}_E H) $,  we obtain the estimate
\beq |T|_h^2(p)  \leq C_{6}, \eeq
where $C_6=LC_4+C_5$.
Since $\mathrm{tr}_EH \leq \rk(E) C$, we get
\beq (|T|_h^2 + L\mathrm{tr}_E H)(p) \leq  C_{7}, \eeq
where $C_7=C_6+\rk(E) LC$.
Therefore, for any point $ x \in M $, one has
\beq |T|^2_h(x) \leq (|T|_h^2 + L\mathrm{tr}_E H)(x)\leq (|T|_h^2 + L\mathrm{tr}_E H)(p) \leq C_{7}. \eeq
By \eqref{c15}, one gets
\beq |\p^{h_0} H|_{h_0}^2 \leq C^2|\p^{h_0}H \cdot H^{-1}  |^2_h \leq  C^2 C_{7}. \eeq
In summary, we establish the uniform $C^1$-estimate
$|H|_{C^1(M,\omega_g, h_0)}\leq C_{8}$ where $ C_{8}$  depends on $M,\omega_g,h_0,  C, c_1$ and $c_2$.
\eproof

\vskip 2\baselineskip

\section{Proofs of Theorem \ref{main2} and Theorem \ref{main4}}
In this section, we prove Theorem \ref{main2} and Theorem \ref{main4}. The weak compactness established in Theorem~\ref{main4} will play a key role in the sequel:

\btheorem
\label{C^0toC^infty} Let $ (M,\omega_g) $ be a compact Hermitian manifold and $ (E, h_0)$ be a Hermitian vector bundle over $M$. Fix $\lambda\in \R$.  Suppose that the following conditions are satisfied  \bd \item
$\{ P_m \} \subset \mathrm{Herm}(E)$
converges smoothly to some $P \in \mathrm{Herm}(E)$;
\item $\{\lambda_m\}\subset \R$ converges to $\lambda$;
\item  $\{h_m\} \subset \mathrm{Herm}^+(E)$ satisfies
\beq \label{curvatureHm}
\Lambda_{\omega_g} \left(\sqrt{-1} R^{h_m}\right) = \lambda_m h_m + P_m;\eeq
\item  there exists a uniform constant $C>0$ such that
\begin{equation}\label{hm-bdd}
C^{-1} h_0 \leq h_m \leq C h_0,
\end{equation}
\ed  holds for all $m$.
Then, there exists a subsequence $\{h_{m_i}\}$ which converges smoothly to a smooth Hermitian metric $h$ on $E$, and this limit satisfies
\begin{equation}\label{curvatureHlimit}
\Lambda_{\omega_g}\left(\sq  R^{h} \right)= \lambda h + P.
\end{equation}
\etheorem

\bproof  Without loss of generality, we assume that 
\beq |\lambda_m|\leq c_1 \eeq 
for all $m$ where $c_1=|\lambda|+1$.
We fix several notations in $\Gamma(M,E^*\otimes E) $:
\beq H_m = h \cdot h_0^{-1}, \quad \Phi^{\lambda_m}_m = \left( \Lambda_{\omega_g} \left(\sqrt{-1} R^{h_m}\right) -\lambda_m h_m\right)\cdot h_0^{-1}, \quad  P^{h_0} = P \cdot h_0^{-1}. \eeq
In this setting, we have 
\beq \Phi_m^{\lambda_m}=P_m\cdot h_0^{-1}.\eeq 
 Since $\lim_{m}\|P_m-P\|_{C^\infty(M, \omega_g,h_0)}=0 $, one has \beq \lim_m\|\Phi^{\lambda_m}_m-P^{h_0}\|_{C^\infty(M, \omega_g,h_0)}=0. \eeq
  In particular,  there exist uniform constants $ c_2> 0 $ and $m_0>0$, depending on $M,\omega_g, h_0$ and $P$, such that for all $ m\geq m_0 $, the following estimate holds
\beq  \nm{\Phi^{\lambda_m}_m}_{C^1(M,\omega_g, h_0)} \leq c_2. \label{close1} \eeq
By Theorem \ref{C1estimate}, there exists $ C_1 = C_1(M,\omega_g,h_0,C, c_1, c_2)$ such that
\beq \nm{H_m}_{C^1(M,\omega_g, h_0)} \leq C_1. \label{close2}\eeq
By \eqref{curvaturedifference},
\be  \Lambda_{\omega_g}\left(\smo\Theta^{h_m}\right) \cdot H_m =\Lambda_{\omega_g}\smo \bp\left(\p^{h_0}H_m \cdot H_m^{-1}\right) \cdot H_m + \Lambda_{\omega_g}\left(\smo\Theta^{h_0}\right) \cdot H_m. \ee
By the definition of $\Phi_m^{\lambda_m}$, one has
\beq \Phi^{\lambda_m}_m + \lambda_m H_m = \Lambda_{\omega_g}\left(\smo\Theta^{h_m}\right) \cdot H_m.\eeq 
If we set $\Om^{\lambda_m}=\Lambda_{\omega_g}\left(\smo\Theta^{h_0}\right)-\lambda_m \mathrm{Id}_E$, then
\beq \Phi^{\lambda_m}_m= \Omega^{\lambda_m}\cdot H_m+ \Lambda_{\omega_g}\smo \bp\left(\p^{h_0}H_m \cdot H_m^{-1}\right) \cdot H_m,\label{close3} \eeq
and so
\beq \label{ellipticHn} \Delta_{\p^{h_0}} H_m =\Phi^{\lambda_m}_m   - \Omega^{\lambda_m} \cdot H_m - \smo\Lambda_{\omega_g}\left(\p^{h_0} H_m \cdot H_m^{-1} \cdot \bar\p H_m\right) . \eeq
Define \beq  W^{\lambda_m}_m:=\Phi^{\lambda_m}_m   - \Omega^{\lambda_m} \cdot H_m - \smo\Lambda_{\omega_g}\left(\p^{h_0} H_m \cdot H_m^{-1} \cdot \bar\p H_m\right).\label{W_mdef}\eeq
Hence, \beq \Delta_{\p^{h_0}} H_m=W^{\lambda_m}_m.\eeq  By  \eqref{close2}, we deduce that 
\beq \left\|\Om^{\lambda_m}\cdot H_m\right\|_{C^0(M,\omega_g, h_0)}\leq C_2\label{close6}\eeq
where  $ C_2 = C_2(M,\omega_g,h_0,C,c_1, c_2)$. By \eqref{hm-bdd}, \eqref{close1}, \eqref{close2} and \eqref{close6}, one deduces that  there exists $ C_3 = C_3(M,\omega_g,h_0,C,c_1, c_2)$  such that
\beq \|W^{\lambda_m}_m\|_{C^0(M,\omega_g, h_0)}\leq C_3. \eeq
 For a fixed large integer $ p $, there exists $ C_4 = C_4(M,\omega_g,h_0,C,c_1, c_2,p)  $ such that 
\beq \nm{W^{\lambda_m}_m}_{L^p(M,\omega_g, h_0)} \leq C_4. \eeq 
By the $W^{2,p}$-estimate for elliptic equations, there exists $ C_5 = C_5(M,\omega_g,h_0,C,c_1, c_2,p)  $ such that
\beq \nm{H_m}_{W^{2,p}(M,\omega_g, h_0)} \leq C_5\left(\nm{H_m}_{L^p(M,\omega_g,h_0)} + \nm{W^{\lambda_m}_m}_{L^p(M,\omega_g,h_0)}\right). \eeq 
In particular, $ \{ H_m \} $ is uniformly bounded in $ W^{2,p}(M, \omega_g, h_0)$.
Since $ p $ is large enough and $M$ is compact, by the compact embedding theorem, there exist a subsequence $\{H_{m_i}\}$, $\alpha\in (0,1)$ and $H\in C^{1,\alpha}(M, \omega_g, h_0)$ such that  \beq \lim_i\| H_{m_i}-H\|_{C^{1,\alpha}(M, \omega_g, h_0)}=0.\label{Hc1alpha}\eeq Moreover,  $ H$ is $h_0$-Hermitian and positive-definite.  Since $\lambda_m\>\lambda$ and 
\beq \Phi_m^{\lambda_m}=P_m\cdot h_0^{-1}\>P^{h_0}=P\cdot h_0^{-1}\eeq 
in $C^\infty$, 
 by using \eqref{W_mdef},  one has \beq  \lim_i \|W^{\lambda_{m_i}}_{m_i}- W\|_{C^{0,\alpha}(M, \omega_g, h_0)}=0\eeq  where $ W \in C^{0,\alpha}(M, \omega_g, h_0) $ is given by 
\beq  W=P^{h_0} - \Omega^\lambda \cdot H - \smo\Lambda_{\omega_g}(\p^{h_0} H \cdot H^{-1} \cdot \bar\p H),\label{Wdef}\eeq
and
\beq \Om^{\lambda}:=\Lambda_{\omega_g}\left(\smo\Theta^{h_0}\right)-\lambda \mathrm{Id}_E. \eeq 
   By Schauder's estimate and  the standard bootstrap process, one concludes that $\{H_m\}$ converges to $H$ in $C^\infty$.  Letting $m_i\>+\infty$  in \eqref{close3}, one has
\beq P^{h_0}= \Omega^\lambda \cdot H + \smo\Lambda_{\omega_g}\bp(\p^{h_0}H \cdot H^{-1}) \cdot H. \label{close4}\eeq
If we set $ h = H \cdot h_0 $, then by \eqref{conformalchange1}, we obtain the Hermitian-Yang-Mills tensor relation:
\beq \Lambda_{\omega_g}\left(\smo\Theta^{h}\right) \cdot H= \Lambda_{\omega_g}\left(\smo\Theta^{h_0}\right)\cdot H + \smo\Lambda_{\omega_g}\bp(\p^{h_0}H \cdot H^{-1}) \cdot H.\label{close5} \eeq
By \eqref{close4} and \eqref{close5}, we  conclude that 
\beq \Lambda_{\omega_g}\left(\smo R^{h}\right) =\lambda h+P. \eeq
This completes the proof. \eproof

\bproof[Proof of Theorem \ref{main2}] We use a similar strategy as in the proof \cite[Theorem~1.1]{WYY26+}.
We define
 a twisted Hermitian-Yang-Mills map
$G : \mathrm{Herm}^+(E) \to \mathrm{Herm}(E)$ by
\begin{equation}\label{eq:Gdef}
G(h) := \Lambda_{\omega_g}\left(\sq R^{h}\right) - \lambda h.
\end{equation}
Hence,  Theorem \ref{main2} is equivalent to proving that
\beq G_0: G^{-1}\left(\mathrm{Herm}^+(E)\right)\>\mathrm{Herm}^+(E)\eeq
is a bijection  where $G_0$ is the restriction of  $G$. Thanks to the curvature  condition \eqref{curvaturecondition}:
$$h_0\in G^{-1}\left(\mathrm{Herm}^+(E)\right).$$
The injectivity of $G_0$ follows from Theorem \ref{comparison}. To
establish the surjectivity, we define \be \mathscr{R}:= \left\{ P
\in \mathrm{Herm}^+(E)\ |\ \text{there exists a  Hermitian metric }
h \text{ and } \Lambda_{\omega_g}\left(\sq R^h\right) - \lambda h =
P \right\}. \ee It is clear that $ \mathscr{R}$ is the image of $G_0$ in
$\mathrm{Herm}^+(E)$.  We show that $\mathscr{R}$ is both relatively
open and  closed in
$\mathrm{Herm}^+(E)$.\\

We establish the openness. For any $ P \in \mathscr{R} $, there exists a Hermitian metric $ h_1 $ on $ E $ such that
\beq \Lambda_{\omega_g}\left(\sq R^{h_1}\right) - \lambda h_1 = P. \eeq
One can show that the tangent map
\beq TG_{h_1}:T_{h_1}\mathrm{Herm}^+(E) \> T_{P}\mathrm{Herm}(E) \eeq
is a bijection. Indeed, we consider the conjugate map $ G_1: \mathrm{Herm}^+(E,h_1) \> \mathrm{Herm}(E,h_1) $ of $G$:
\beq G_1(H) = G(H \cdot h_1) \cdot h_1^{-1}. \eeq
It is clear that
\beq G_1(H) = \left(\Omega_1 +  \Lambda_{\omega_g}\smo \bp\left(\p^{h_1}H \cdot H^{-1}\right) \right) \cdot H, \eeq
where $ \Omega_1 \in \Gamma(M,E^*\otimes E) $ is given by
\beq \Omega_1 :=  P \cdot h_1^{-1} . \eeq
Therefore, the linearization of $G_1$ at point $ \mathrm{Id}_E $ is
\beq \mathscr{L}_{\mathrm{Id}_E}(\Psi) = \Delta_{\p^{h_1}} \Psi + \Omega_1 \cdot \Psi. \eeq
By the elliptic theory described in the proof of \cite[Theorem~1.1]{FWYY26+}, $\mathscr{L}_{\mathrm{Id}_E}$ is bijective. Consequently, $TG_{h_1}$ is also bijective, and the openness follows from the implicit function theorem.\\

 The closeness follows from  Theorem \ref{C^0toC^infty}. Indeed, let $ \{ P_m \}\subset \mathscr{R} $ be a sequence such that
$\{P_m\}$ converges to   $ P \in \mathrm{Herm}^+(E) $ in  $C^{\infty}(M,E^*\otimes E) $. In particular,  there exist Hermitian metrics $ h_m $ such that
 \beq \Lambda_{\omega_g}\left(\sq R^{h_m}\right) = \lambda h_m + P_m. \eeq
 Moreover, there exists a constant $C=C(\omega_g, P, h_0)$ such that
 \beq C^{-1}P_0 \leq P_m \leq C P_0,  \eeq
for large $m$.  That is,
\beq C^{-1}\left(\Lambda_{\omega_g}\left(\sq R^{h_0}\right) - \lambda h_0\right)\leq \Lambda_{\omega_g}\left(\sq R^{h_m}\right) - \lambda h_m\leq C\left(\Lambda_{\omega_g}\left(\sq R^{h_0}\right) - \lambda h_0\right).\eeq
By Theorem \ref{comparison},
 \beq C^{-1}h_0 \leq h_m \leq C h_0. \eeq
  By Theorem \ref{C^0toC^infty}, there exist a subsequence $ \{h_{m_i}\}$ and a smooth Hermitian metric $ h $  on $E$ such that $ \{h_{m_i}\} $ converges to $ h $ in $ C^\infty(M,E^*\otimes E) $ and  \beq \Lambda_{\omega_g}\left(\sq R^h\right) = \lambda h + P. \eeq
Hence, $ P \in \mathscr{R} $ and so $ \mathscr{R} $ is relatively closed in $\mathrm{Herm}^+(E) $. This completes the proof.
\eproof

\vskip 2\baselineskip

\section{Proof of Theorem \ref{main}}

In this section we prove Theorem \ref{main}:

  \btheorem
Let $ (M,\omega_g) $ be a compact Gauduchon manifold and $ E $  be an arbitrary holomorphic vector bundle over $ M $.  If we define 
\beq \mu_E^-: = \sup
\left\{\ \lambda\ \left|\  \text{ there exists a Hermitian metric }
h \text{ such that } \Lambda_{\omega_g}\left(\sq R^h\right) >
\lambda h \right.\right\}, \label{mu} \eeq
then \beq \lambda_E^-=\mu_E^-.\eeq  Moreover,   for every $\lambda<\lambda_E^-$ and every positive-definite Hermitian tensor $P\in \Gamma(M,E^*\ts \bar E^* )$, there exists a unique Hermitian metric $h$ on $E$ satisfying the twisted Hermitian-Einstein equation
\beq \Lambda_{\omega_g}\left(\sq R^h\right)=\lambda\cdot h+ P.\eeq
\etheorem

\bproof Recall that, for a coherent quotient sheaf $\sQ$ of $E$,
\beq \lambda_\sQ=\frac{2\pi n\int_M c^{\mathrm{BC}}_1(\sQ)\wedge
	\omega_g^{n-1}}{\mathrm{\rk}(\sQ)\int_M\omega_g^n}=\frac{\deg_{\omega_g}(\sQ)}{\rk(\sQ)
	\int_M\omega_g^n}, \eeq  and \beq \lambda_E^- := \inf\left\{ \left.
\lambda_\sQ\ \right|\  \mathscr{Q} \text{ is a coherent quotient
	sheaf of } E, \mathrm{rk}(\mathscr{Q}) > 0 \right\}. \eeq  Since $M$
is compact, it is easy to see that $\mu_E^-$ is well-defined and  $
\mu_E^- \leq \lambda_E^-$.  Suppose for contradiction that
$\mu_E^->\lambda_E^-$.  There exists some $\lambda>\lambda_E^-$ such
that \beq \Lambda_{\omega_g}\left(\sq R^h\right)\geq \lambda h.
\label{cur}\eeq
Let $\sQ$ be a coherent quotient sheaf of $E$ with $\rk(\sQ)>0$ and
\beq \lambda_E^-\leq \lambda_\sQ<\lambda.\eeq Since the curvature
increases along quotients, by \eqref{cur},  one gets \beq
\deg_{\omega_g}(\sQ)\geq \lambda \cdot  \rk(\sQ) \int_M \omega_g^n.
\eeq In particular, \beq
\lambda_\sQ=\frac{\deg_{\omega_g}(\sQ)}{\rk(\sQ)
	\int_M\omega_g^n}\geq \lambda.\eeq
This is a contradiction.\\

We claim that
\beq\label{Claimlambda} \mu_E^- = \lambda_E^-. \eeq Suppose, for
contradiction, that  $ \mu_E^- < \lambda_E^- $. By definition of $
\mu_E^- $, there exist sequences $ \{ \mu_m \}_{m=0,1,\cdots}
\subset \mathbb{R} $ and $ \{ g_m \} \subset \mathrm{Herm}^+(E) $
such that \beq \mu_m\nearrow\mu_E^-, \eeq  and \beq
\Lambda_{\omega_g}\left(\sq R^{g_m}\right) > \mu_m g_m. \eeq On the
other hand, by Theorem \ref{main2},  for any $ P \in
\mathrm{Herm}^+(E) $, there exist Hermitian metrics $ \{h_m\}$ on
$E$ such that \beq \label{lambdamsol} \Lambda_{\omega_g}\left(\sq
R^{h_m}\right) = \mu_m h_m+P. \eeq We choose  $ h_0 $ to be the
reference Hermitian metric on $E$. In particular,
$$\Lambda_{\omega_g}\left(\sq R^{h_0}\right) = \mu_0 h_0+P.$$
Moreover, for any $ m \in \mathbb{N} $, one has
\beq \Lambda_{\omega_g}\left(\sq R^{h_m}\right)  - \mu_0 h_m =
(\mu_m - \mu_0)h_m + P \geq  P>0, \eeq since the sequence
$\{\mu_m\}$ is increasing. Therefore,  for any $ m \in \mathbb{N} $,
\beq \Lambda_{\omega_g}\left(\sq R^{h_0}\right)  - \mu_0 h_0 \leq
\Lambda_{\omega_g}\left(\sq R^{h_m}\right)  - \mu_0 h_m, \quad
\Lambda_{\omega_g}\left(\sq R^{h_m}\right)  - \mu_0 h_m > 0. \eeq By
Theorem \ref{comparison}, the following inequality holds \beq h_0
\leq h_m \label{lowerbound}\eeq
for all $m$. \\

We discuss the uniform upper bound of $h_m$ and define \beq
H_m=h_m\cdot h_0^{-1}\in\Gamma(M,E^*\ts E). \eeq Let  \beq \Lambda_m
:= \sup_M \lambda_{\max}(H_m)\eeq denote the supremum over $M$ of
the largest eigenvalue of $H_m$.  Define \beq \tilde h_m:=
\Lambda_m^{-1}\cdot h_m, \quad \tilde H_m:=\Lambda_m^{-1}\cdot H_m.
\eeq
Then it is clear that $ \tilde H_m\leq \mathrm{Id}_E$ with respect  to both $h_0$ and $h_m$.\\

$(1)$. Suppose that \beq \limsup_{m\>\infty}\Lambda_m = \limsup_{m\>\infty}\sup_M \lambda_{\max}(H_m)=+\infty.\eeq
By \eqref{lambdamsol}, we  deduce that \beq
\label{lambdasoltildeh}\Lambda_{\omega_g}\left(\sq \Theta^{\tilde
	h_m}\right) =\Lambda_{\omega_g}\left(\sq \Theta^{h_m}\right) = \mu_m
\mathrm{Id}_E+P \cdot h_m^{-1} \eeq as sections in $
\Gamma(M,E^*\otimes E) $. Since $h_0\leq h_m$, there exists a
constant $ C_1 = C_1(P, h_0) > 0 $ such that for any $ m \in
\mathbb{N} $, \beq P\leq C_1 h_0\leq C_1 h_m.\eeq In particular,
\beq \mu_0 \mathrm{Id}_E \leq \Lambda_{\omega_g}\left(\sq
\Theta^{\tilde h_m}\right)  \leq
\left(C_1+\mu_E^-\right)\mathrm{Id}_E, \label{boundedcurvature} \eeq
and $0<\tilde H_m\leq \mathrm{Id}_E$
as $ \tilde h_m $-Hermitian sections.

For any $\sigma\in (0,1]$, one has the following inequality: \beq
\left| \p^{h_0}H_m^\sigma\cdot H_m^{-\sigma/2} \right|_{g,h_0}^2
\leq  \frac{1}{\sigma}\Delta_\C\left(\tr_E H_m^\sigma\right)
+\tr_E\left(\left(\Lambda_{\om_g}
\left(\sq\Theta^{h_m}\right)-\Lambda_{\om_g}\left(\sq
\Theta^{h_0}\right)\right)\cdot H_m^\sigma\right). \label{estimateA}
\eeq Indeed, it follows from \eqref{conformalchange1} and an
elementary calculation that \beq
\Lambda_{\om_g}\sq\left\{(\p^{h_0}H)H^{-1},\p^{h_0}H^\sigma\right\}_{h_0}\geq
\left| \p^{h_0}H_m^\sigma\cdot H_m^{-\sigma/2}
\right|_{g,h_0}^2.\eeq Moreover, by \eqref{lambdamsol} \beq
\Lambda_{\om_g}
\left(\sq\Theta^{h_m}\right)-\Lambda_{\om_g}\left(\sq
\Theta^{h_0}\right)=(\mu_m-\mu_0) \mathrm{Id}_E+P\cdot
(h_m^{-1}-h_0^{-1}). \eeq Since $h_0\leq h_m$, we deduce that
\beq\left|\tr_E\left(\left(\Lambda_{\om_g}
\left(\sq\Theta^{h_m}\right)-\Lambda_{\om_g}\left(\sq
\Theta^{h_0}\right)\right)\cdot H_m^\sigma\right)\right|\leq
C_2\tr_EH_m^\sigma, \eeq where $C_2=C_2(M, \omega_g, P,\mu_E^-,
\mu_0)$.
Hence, we get
\beq \left|\p^{h_0}H_m^\sigma \cdot H_m^{-\sigma/2}\right|_{g,h_0}^2
\leq  \frac{1}{\sigma}\Delta_\C\left(\tr_E H_m^\sigma\right)
+C_2\tr_E H_m^\sigma. \eeq Multiplying by $\Lambda_m^{-\sigma}$
yields \beq \left|\p^{h_0}\tilde H_m^\sigma \cdot \tilde
H_m^{-\sigma/2}\right|_{g,h_0}^2 \leq
\frac{1}{\sigma}\Delta_\C\left(\tr_E \tilde H_m^\sigma\right)
+C_2\tr_E \tilde H_m^\sigma. \label{estimateHsigma} \eeq Since
$0<\tilde H_m\leq\id_E$, we have $\tr_E\tilde H_m^\sigma\leq\rk(E)$
and $\tilde H_m^{-\sigma/2}\geq\id_E$. This yields \beq
\left|\p^{h_0} \tilde H_m^\sigma \right|_{g,h_0}^2 \leq
\left|\p^{h_0}\tilde H_m^\sigma \cdot \tilde
H_m^{-\sigma/2}\right|_{g,h_0}^2 \leq
\frac{1}{\sigma}\Delta_\C\left(\tr_E \tilde H_m^\sigma\right) +C_3.
\label{estimateHsigma1} \eeq  with $C_3=C_2\rk(E)$. Since $0<\tilde
H_m\leq\id_E$, \beq \left\| \tilde H_m^\sigma
\right\|_{L^2(M,\omega_g,h_0)}^2 \leq \rk(E)\mathrm{Vol}(M,\om_g).
\label{tilde H L2 estimate} \eeq Integrating inequality
\eqref{estimateHsigma1} over $M$ and using the Gauduchon condition
$\partial\bar{\partial}\,\omega_g^{n-1}=0$, we obtain
\begin{equation}
\label{d tilde H L2 estimate} \bigl\|\partial^{h_0}\tilde
H_m^\sigma\bigr\|_{L^2(M,\omega_g,h_0)}^2 \;\le\;
C_3\,\mathrm{Vol}(M,\omega_g).
\end{equation}
Since $\tilde H_m^\sigma$ is Hermitian with respect to $h_0$, i.e.,
$(\tilde H_m^\sigma)^*=\tilde H_m^\sigma$, we have
\[
\bar{\partial}\tilde H_m^\sigma = \bigl(\partial^{h_0}\tilde
H_m^\sigma\bigr)^*,
\]
and therefore
\begin{equation}
\label{d tilde H equation} \bigl\|\bar{\partial}\tilde
H_m^\sigma\bigr\|_{L^2(M,\omega_g,h_0)} =
\bigl\|\partial^{h_0}\tilde H_m^\sigma\bigr\|_{L^2(M,\omega_g,h_0)}.
\end{equation}
Combining \eqref{tilde H L2 estimate}, \eqref{d tilde H L2 estimate}
and \eqref{d tilde H equation} yields a uniform $W^{1,2}$-bound \beq
\|\tilde H_m^\sigma\|_{W^{1,2}(M,\omega_g,h_0)} \le C_4, \eeq where
$C_4=C_4(M,\omega_g,h_0, P, \mu_E^-)$. By the weak compactness
theorem, there exist a subsequence $\{\tilde H_{m_i}\}$ of $\tilde
H_m$,  and a sequence $\{\tau_i\}\subset  (0,1/2]$ decreasing to $0$
such that

\bd \item [(a)] $\lim_i\Lambda_{m_i}=+\infty$.
\item[(b)] $ \{\tilde H_{m_i}\}$ converges to  some  $\tilde H_\infty$ in the weak $ W^{1,2}(M, \omega_g, h_0) $-sense;
\item[(c)] $ \{\tilde H^{2\tau_i}_\infty\}$ converges to some $ \tilde H $ in the weak $ W^{1,2}(M, \omega_g, h_0) $-sense.
\ed On the other hand,  taking $\sigma=1$  in
\eqref{estimateHsigma}, one has \beq -\Delta_\C \left(\tr_E \tilde
H_m\right)\leq C_2\tr_E\tilde H_m. \eeq Applying the Moser iteration
to this inequality yields \beq \sup_M  \tr_E \tilde H_m \leq
C_5\left\|\tr_E \tilde H_m\right\|_{L^2(M,g,h_0)}, \eeq where
$C_5=C_5(M, \omega_g, C_2)$. By construction, $ \sup_M\tr_E\tilde
H_m\geq1$,  and so \beq 1\leq C_5\left\|\tr_E \tilde
H_m\right\|_{L^2(M,g,h_0)} \leq C_6\left\|\tilde
H_m\right\|_{L^2(M,g,h_0)}, \eeq with $C_6=C_6(M, \om_g,h_0, C_5)$.
Therefore, \beq \left\|\tilde H_m\right\|_{L^2(M,g,h_0)}\geq
C_6^{-1}.\eeq As a consequence,
$$\tilde H\neq 0. $$   By using a construction similar to that in \cite{UY86},
\beq  \pi = \mathrm{Id}_E - \tilde H \eeq  defines a weak
holomorphic projection with respect to  $ h_0 $. That is, \beq
\pi^*=\pi=\pi^2\eeq and \beq  \left(\id_E-\pi\right)\circ \bp\pi =0
\eeq hold almost everywhere on $M$. Moreover, $\pi$ defines a
coherent subsheaf $\cF\subset E$.  Let $\Sigma\subset M$ be the
singular locus of $\cF$. Then
\beq \mathrm{tr}_E\left(\pi|_{M\setminus\Sigma}\right)\equiv \rk(\cF) <\rk(E).\eeq
By computing the second fundamental form, the degree of $\cF$ is
given by: \beq \deg_{\omega_g}(\cF) = \int_{M\setminus \Sigma}
\tr_{E}\left(\Lambda_{\om_g}\left(\sq\Theta^{h_0}\right)\cdot
\pi\right)\om_g^n -\int_{M\setminus \Sigma} \left|B\right|_{g,h_0}^2
\om_g^n, \eeq where $B$ is the second fundamental form of $ \cF\>E$.
A simple calculation shows \beq |B|^2_{g,h_0}=|\p^{h_0}\pi|^2_{g,
	h_0}.\eeq Hence, we obtain \beq \deg_{\omega_g}(\cF) = \int_{M}
\tr_{E}\left(\Lambda_{\om_g}\left(\sq\Theta^{h_0}\right)\cdot
\pi\right)\om_g^n -\int_{M} \left|\partial^{h_0}
\pi\right|_{g,h_0}^2 \om_g^n. \label{degree} \eeq On the other hand,
by a scaling of \eqref{estimateA} and integrating over $M$,  we have
the inequality
\begin{eqnarray}
\int_M  \left|\p^{h_0}\tilde H_{m_i}^{2\tau_j}\right|_{g,h_0}^2
\om_g^n
&\leq& -\int_M\tr_E\left( \Lambda_{\om_g}\left(\sq \Theta^{h_0}\right) \cdot \tilde H_{m_i}^{2\tau_j}\right)\om_g^n\nonumber\\
&&+ \int_M\tr_E\left(\Lambda_{\om_g}\left(\sq
\Theta^{h_{m_i}}\right)\cdot \tilde H_{m_i}^{2\tau_j}\right)\om_g^n.
\end{eqnarray}
By taking limit, one obtains
\begin{eqnarray}\int_M  \left|\p^{h_0}\tilde H\right|_{g,h_0}^2  \om_g^n
&\leq& -\int_M\tr_E\left( \Lambda_{\om_g}\left(\sq \Theta^{h_0}\right) \cdot \tilde H\right)\om_g^n\nonumber\\
&&+ \limsup_{\tau_j \> 0^+}\limsup_{m_i \>
	+\infty}\int_M\tr_E\left(\Lambda_{\om_g}\left(\sq
\Theta^{h_{m_i}}\right)\cdot \tilde
H_{m_i}^{2\tau_j}\right)\om_g^n.\label{limitdegree}
\end{eqnarray}
Note also that $ \p^{h_0}\pi=-\p^{h_0}\tilde H$ holds almost
everywhere on $M$. Hence, by \eqref{degree} and \eqref{limitdegree},
\beq
\deg_{\omega_g}(\cF)
\geq  \deg_{\omega_g}(E) -\limsup_{\tau_j \> 0^+}\limsup_{m_i \>
	+\infty}\int_M\tr_E\left(\Lambda_{\om_g}\left(\sq
\Theta^{h_{m_i}}\right)\cdot \tilde H_{m_i}^{2\tau_j}\right)\om_g^n.
\eeq If we denote $ \mathscr{Q} = E/\mathscr{F} $, then \beq
\label{quotientsheafQ} \mathrm{deg}_{\omega_g}(\mathscr{Q}) \leq
\limsup_{\tau_j \> 0^+}\limsup_{m_i \>
	+\infty}\int_M\tr_E\left(\Lambda_{\om_g}\left(\sq
\Theta^{h_{m_i}}\right)\cdot \tilde H_{m_i}^{2\tau_j}\right)\om_g^n.
\eeq On the other hand,  it is clear that \beq \limsup_{m_i \>
	+\infty} \left(\mu_{m_i}\mathrm{Id}_E , \tilde
H_{m_i}^{2\tau_j}\right)_{h_0}=\mu_E^-\left(\mathrm{Id}_E,\tilde
H^{2\tau_j}_\infty\right)_{h_0},\eeq and \beq \left(P \cdot
h_{m_i}^{-1}, \tilde H_{m_i}^{2\tau_j}\right)_{h_0}=
\Lambda_{m_i}^{-2\tau_j}\left(P \cdot h_{m_i}^{-1},
H_{m_i}^{2\tau_j}\right)_{h_0}. \eeq By \eqref{lowerbound}, we know
$H_{m_i}\geq \mathrm{Id}_E$. Moreover, since $2\tau_j-1\leq 0$, \beq
H_{m_i}^{2\tau_j-1}\leq \mathrm{Id}_E. \eeq This implies \beq
h_0\left(P\cdot h_{m_i}^{-1},H_{m_i}^{2\tau_j}\right)
=h_0\left(P\cdot h_{0}^{-1},H_{m_i}^{2\tau_j-1}\right)\leq
\tr_{h_0}P, \eeq
and
\beq  \limsup_{m_i \> +\infty}\left(P \cdot
h_{m_i}^{-1}, \tilde H_{m_i}^{2\tau_j}\right)_{h_0}=\limsup_{m_i \> +\infty} \Lambda_{m_i}^{-2\tau_j}\left(P \cdot h_{m_i}^{-1},  H_{m_i}^{2\tau_j}\right)_{h_0}=0.\eeq
By \eqref{lambdasoltildeh},
\be  \limsup_{m_i \> +\infty} \left(\Lambda_{\omega_g}\left(\sq \Theta^{\tilde h_{m_i}}\right), \tilde H_{m_i}^{2\tau_j}\right)_{h_0} &=& \limsup_{m_i \> +\infty} \left(\mu_{m_i}\mathrm{Id}_E + P \cdot h_{m_i}^{-1}, \tilde H_{m_i}^{2\tau_j}\right)_{h_0} \\
& \leq & \limsup_{m_i \> +\infty} \left(\mu_{m_i}\mathrm{Id}_E, \tilde H_{m_i}^{2\tau_j}\right)_{h_0} + \limsup_{m_i \> +\infty} \left(P \cdot h_{m_i}^{-1}, \tilde H_{m_i}^{2\tau_j}\right)_{h_0} \\
& = & \mu_E^-\left(\mathrm{Id}_E,\tilde H^{2\tau_j}_\infty\right)_{h_0} .
\ee
When $ \tau_j \> 0^+ $, one  deduces that \be\limsup_{\tau_j\>0+}
\limsup_{m_i \> +\infty} \left(\Lambda_{\omega_g}\left(\sq
\Theta^{\tilde h_{m_i}}\right), \tilde
H_{m_i}^{2\tau_j}\right)_{h_0}
& \leq & \limsup_{\tau_j \> 0^+}\mu_E^-\left(\mathrm{Id}_E,\tilde H^{2\tau_j}_\infty\right)_{h_0} \\
& = & \mu_E^-\left(\mathrm{Id}_E, \tilde H\right)_{h_0} =
\mu_E^-\mathrm{rk}(\mathscr{Q})\left(\int_M \omega_g^n\right). \ee
Hence, by \eqref{quotientsheafQ}, one obtains \beq
\label{slopeQestimate}  \mathrm{deg}_{\omega_g}(\mathscr{Q}) \leq
\mu_E^-\mathrm{rk}(\mathscr{Q})\left(\int_M \omega_g^n\right), \eeq
By \eqref{slopeQestimate} and the definition of $ \lambda_E^- $,  we
obtain \beq \lambda_E^-\leq \lambda_\sQ\leq \mu_E^-. \eeq
This contradicts the assumption that $ \mu_E^- < \lambda_E^- $.\\

$(2)$. Suppose that
there exists a uniform constant $ C > 0 $ such that
\beq \limsup_{m\>\infty}\Lambda_m = \limsup_{m\>\infty}\sup_M
\lambda_{\max}(H_m)=C<+\infty.\eeq Then we have \beq h_m \leq (C+1)
h_0, \eeq for all large $ m$.  Moreover, by \eqref{lowerbound},
there exists some $C_7>0$ such that \beq C_7^{-1} h_0\leq  h_m \leq
C_7 h_0. \eeq By \eqref{lambdamsol}, we have \beq
\Lambda_{\omega_g}\left(\sq R^{h_m}\right) = \mu_m h_m + P. \eeq By
Theorem \ref{C^0toC^infty}, there exist a subsequence $ \{
h_{m_j}\}$ and a Hermitian metric $ h $ on $ E $ such that
$\{h_{m_j}\}$ converges to  $h$ in $ C^\infty(M,E^*\otimes E)$ and
\beq \Lambda_{\omega_g}\left(\sq R^h\right) = \mu_E^- h + P. \eeq
Since $ P > 0 $, there exists $ \eps > 0 $ such that
\beq\Lambda_{\omega_g}\left(\sq R^h\right) - (\mu_E^-+\eps)h =
P-\eps h
> 0. \eeq
This contradicts to the definition of $ \mu_E^- $. \\

Hence, we deduce that $ \mu_E^- = \lambda_E^- $. In particular, for
any $ \lambda < \lambda_E^- $, there exists a $ \lambda_1 \in
(\lambda, \lambda_E^-) $ and a Hermitian metric $ h_1 $ on $ E $
such that \beq \Lambda_{\omega_g}\left(\sq R^{h_1}\right)> \lambda_1
h_1 > \lambda h_1. \eeq By Theorem \ref{main2}, for any $ P \in
\mathrm{Herm}^+(E) $, there exists a unique Hermitian metric $ h $
on $ E $ such that \beq \Lambda_{\omega_g}\left(\sq R^{h}\right)=
\lambda h + P. \eeq This completes the proof. \eproof

\vskip 2\baselineskip

\section{The Chern number inequality and proofs of Theorems \ref{main7} and \ref{main6}}
In this section, we prove Theorem \ref{main7} and Theorem \ref{main6}. Let $E$ be a holomorphic vector bundle over a compact Gauduchon manifold $(M,\omega_g)$.  We define 
 \beq \mu_E^+: = \inf
\left\{\lambda\left|\  \text{ there exists a Hermitian metric }
h \text{ such that } \Lambda_{\omega_g}\left(\sq R^h\right) <
\lambda h \right.\right\}.\label{mu+} \eeq

\blemma\label{dual}  Let $ (M,\omega_g) $ be a compact Gauduchon manifold and $ E $ be a holomorphic vector bundle over $ M $. Then 
\beq \mu_E^+=-\mu_{E^*}^-=-\lambda_{E^*}^-=\lambda_E^+. \eeq 
\elemma 

\bproof  By duality,  coherent subsheaves of $E$   correspond to quotient sheaves of $E^*$.  Hence, 
\beq \lambda_E^+=-\lambda_{E^*}^-.\label{dual2}\eeq 
Moreover, for any $\lambda>\mu_E^+$, 
there exists a Hermitian metric $h^E$ on $E$  such that 
\beq \Lambda_{\omega_g}\left(\sq R^{h^E}\right) <
\lambda h^E .\eeq 
This is equivalent to 
\beq \Lambda_{\omega_g}\left(\sq \Theta^{h^E}\right) <
\lambda \mathrm{Id}_E\eeq
with respect to $h^E$. By duality, the dual metric $h^{E^*}$ on $E^*$ satisfies  
\beq \Lambda_{\omega_g}\left(\sq \Theta^{h^{E^*}}\right) =-\Lambda_{\omega_g}\left(\sq \Theta^{h^E}\right)^t \in \Gamma(M,E\ts E^*).\eeq
Hence, 
\beq \Lambda_{\omega_g}\left(\sq \Theta^{h^{E^*}}\right) >-\lambda \mathrm{Id}_{E^*}.\eeq
By the definition of $\mu_{E^*}^-$, one obtains
\beq -\lambda \leq \mu_{E^*}^-. \eeq 
In particular, letting $\lambda\>\mu_E^+$, we have
\beq  -\mu_E^+ \leq \mu_{E^*}^-.\eeq 
Similarly, for any $\lambda<\mu_{E^*}^-$, 
there exists a Hermitian metric $h^{E^*}$ on $E^*$  such that 
\beq \Lambda_{\omega_g}\left(\sq R^{h^{E^*}}\right) >
\lambda h^{E^*} .\eeq
By duality, one has  
\beq \Lambda_{\omega_g}\left(\sq R^{h^E}\right) <-
\lambda h^E .\eeq 
Hence, $-\lambda\geq \mu_{E}^+$ and so $-\mu_{E^*}^- \geq \mu_{E}^+$. Hence, $\mu_{E}^+=-\mu_{E^*}^-$. Moreover, by Theorem \ref{main}, $\lambda_{E^*}^-=\mu_{E^*}^-$. By \eqref{dual2}, one gets the conclusion.
\eproof 

\noindent 
The following is Theorem \ref{main7}:
\btheorem
\label{epspinch}
Let $ (M,\omega_g) $ be a compact Gauduchon manifold and $ E $ be a holomorphic vector bundle over $ M $. For any $ \eps > 0 $, there exists a Hermitian metric $ h_\eps $ on $ E $ such that 
\beq \left(\lambda_E^-- \eps\right) \cdot h_\eps \leq \Lambda_{\omega_g}\left(\smo R^{h_\eps}\right) \leq \left(\lambda_E^++ \eps\right) \cdot h_\eps. \eeq

\etheorem

\bproof

Fix $ \eps > 0 $. By Lemma \ref{dual}, $\mu_E^+=\lambda_E^+$, and by the definition of $\mu_E^+$, there exists a Hermitian metric $ g_\eps $ on $ E$ such that 
\beq \Lambda_{\omega_g}\left(\smo R^{g_\eps}\right) < (\lambda_E^++ \eps) g_\eps. \eeq
In particular, 
\beq \label{Identitygeps} \Lambda_{\omega_g}\left(\smo R^{g_\eps}\right) - (\lambda_E^- - \eps)g_\eps < (\lambda_E^+- \lambda_E^-+ 2\eps) g_\eps. \eeq 
Note that $\lambda_E^+\geq  \lambda_E^-$. 
On the other hand, by Theorem \ref{main}, there exists a Hermitian metirc $ h_\eps $ on $E$ such that 
\beq \label{gheps} \Lambda_{\omega_g}\left(\smo R^{h_\eps}\right) = \left(\lambda_E^- - \eps\right)h_\eps + g_\eps. \eeq
In particular, 
\beq \label{Identityheps}\Lambda_{\omega_g}\left(\smo R^{h_\eps}\right) - \left(\lambda_E^--\eps\right)h_\eps = g_\eps > 0. \eeq
By \eqref{Identitygeps} and \eqref{Identityheps}, one has 
$$ \Lambda_{\omega_g}\left(\smo R^{g_\eps}\right) - \left(\lambda_E^- -\eps\right)g_\eps \leq (\lambda_E^+- \lambda_E^-+ 2\eps)\left(\Lambda_{\omega_g}\left(\smo R^{h_\eps}\right) - \left(\lambda_E^--\eps\right)h_\eps\right). $$
Thanks to Theorem \ref{main5}, it follows that 
\beq \label{comparisongheps} 0<g_\eps \leq (\lambda_E^+- \lambda_E^-+ 2\eps)h_\eps. \eeq
By using the estimate \eqref{comparisongheps}  for $g_\eps$ in  \eqref{gheps}, we conclude that 
\beq \left(\lambda_E^- - \eps\right) \cdot h_\eps \leq \Lambda_{\omega_g}\left(\smo R^{h_\eps}\right) \leq \left(\lambda_E^++ \eps\right) \cdot h_\eps. \eeq
This completes the proof.
\eproof

\bproof[Proof of Theorem \ref{main6}] 	
Let $\{z^i\}$ be local holomorphic coordinates on $M$ and $\{e_\alpha\}$ be holomorphic frames of $E$.  Let $R^h\in \Gamma(M,\Lambda^{1,1}T^*M\ts E^*\ts \bar E^*)$ be the Chern curvature  tensor of $(E,h)$.
We set $$R^{(1)}_{i\bar j} =h^{\alpha\bar\beta}R^h_{i\bar j \alpha\bar \beta}, \quad R^{(2)}_{\alpha\bar\beta}=g^{i\bar j}  R^h_{i\bar j \alpha\bar\beta } \qtq{and} s_h =g^{i\bar j} R^{(1)}_{i\bar j}.$$
The following expressions are essentially  well-known:
\beq \label{c1computation} \int_M c^2_1(E) \wedge \omega_g^{n-2} = \frac{1}{4\pi^2n(n-1)}\int_M \left(s_h^2 - \left|\mathrm{Ric}^{(1)}\right|_g^2 \right) \omega_g^n, \eeq
and 
\beq \label{c2computation} \int_M c_2(E) \wedge \omega_g^{n-2} = \frac{1}{8\pi^2n(n-1)}\int_M \left(s_h^2 - \left|\mathrm{Ric}^{(1)}\right|_g^2 - \left| \mathrm{Ric}^{(2)}\right|^2_h + |R^h|_{g,h}^2\right) \omega_g^n. \eeq
We define $T\in \Gamma(M,\Lambda^{1,1}T^*M\ts E^*\ts \bar E^*)$:
\beq T_{i\bar j \alpha\bar\beta } = R^h_{i\bar j\alpha\bar\beta} - \frac{1}{n}g_{i\bar j}R^{(2)}_{\alpha\bar\beta} - \frac{1}{r}R^{(1)}_{i\bar j}h_{\alpha \bar\beta} + \frac{1}{nr}g_{i\bar j}h_{\alpha \bar\beta}s_h. \eeq
A straightforward calculation shows that 
\beq |T|^2=|R^h|^2 - \frac{1}{n}\left|\mathrm{Ric}^{(2)}\right|^2 - \frac{1}{r}\left|\mathrm{Ric}^{(1)}\right|^2 + \frac{1}{nr} s_h^2.\eeq 
One has
\be && \int_M \left(2rc_2(E) - (r-1)c^2_1(E)\right) \wedge \omega_g^{n-2} \\
& = & \frac{1}{4\pi^2n(n-1)}\int_M \left(s_h^2 - \left|\mathrm{Ric}^{(1)}\right|^2 - r\left|\mathrm{Ric}^{(2)}\right|^2 + r|R^h|^2\right) \omega_g^n \\
& = & \frac{1}{4\pi^2n(n-1)}\int_M \left(r|T|^2 + \frac{n-1}{n}\left( s_h^2 - r\left|\mathrm{Ric}^{(2)}\right|^2\right) \right) \omega_g^n, \ee
where we substitute  $|T|^2$ in the second identity.
Let $\{\lambda_1,\cdots, \lambda_r\}$ be the eigenvalues of $ \Lambda_{\omega_g}\left(\sq R^h\right) $  with respect to $ h $.  If we assume that  $$ a h\leq \Lambda_{\omega_g}\left(\sq R^h\right) \leq bh, $$ for $a,b \in \R$, then
$ a \leq \lambda_k \leq b$. 
A straightforward combinatorial  argument shows that
\beq  \sum_{i<j} (\lambda_i-\lambda_j)^2 \leq \Bigl\lfloor \frac{r^2}{4} \Bigr\rfloor(b-a)^2.\eeq 
This yields
\beq s_h^2 - r\left|\mathrm{Ric}^{(2)}\right|^2 = \left(\sum_{k=1}^r \lambda_k\right)^2 - r\sum_{k=1}^r \lambda_k^2=- \sum_{i<j} (\lambda_i-\lambda_j)^2 \geq- \Bigl\lfloor \frac{r^2}{4} \Bigr\rfloor(b-a)^2. \eeq
Hence, the following Chern number inequality holds:
$$ \int_M \left(2rc_2(E) - (r-1)c^2_1(E)\right) \wedge \omega_g^{n-2}\geq  -  \Bigl\lfloor \frac{r^2}{4} \Bigr\rfloor\frac{\left(b - a\right)^2}{4\pi^2n^2}\int_M \omega_g^n+\frac{r}{4\pi^2n(n-1)}\int_M |T|^2\omega_g^n.$$
By Theorem \ref{epspinch}, for any $\eps>0$, we have $a=\lambda_E^--\eps$, $b=\lambda_E^++\eps$ and $h=h_\eps$. Hence, 
$$ \int_M \left(2rc_2(E) - (r-1)c^2_1(E)\right) \wedge \omega_g^{n-2}\geq  -  \Bigl\lfloor \frac{r^2}{4} \Bigr\rfloor\frac{\left(\lambda_E^+-\lambda_E^- +2\eps\right)^2}{4\pi^2n^2}\int_M \omega_g^n.$$
Letting $\eps\>0^+$, the conclusion follows.
\eproof

\vskip 2\baselineskip
\section{ Proof of Theorem \ref{main3}}

In this section, we  provide an alternative proof of the Donaldson-Uhlenbeck-Yau theorem for stable vector bundles and prove Theorem \ref{main3}:

\btheorem Let $ (M,\omega_g) $ be a compact  Gauduchon manifold and
$ E $ be a  holomorphic vector bundle over $ M $. If $E$ is
semi-stable, then for any $ P \in \mathrm{Herm}^+(E) $ and for any $
\eps > 0 $, there exists a unique smooth Hermitian metric $ h_\eps $
on $E$ such that \beq \Lambda_{\omega_g}\left(\smo R^{h_\eps}\right)
= (\lambda_E - \eps)h_\eps + \eps P, \eeq and $ E $ is approximately
Hermitian-Einstein. Moreover, if $E$ is stable, then a subsequence
$\{h_{\eps_i}\}$ converges smoothly to a Hermitian-Einstein metric
$h$ on $E$ as $\eps_i\>0$. \etheorem

\bproof It is clear that for any vector bundle $ \lambda_E\geq
\lambda_E^-$. Since $ E $ is semi-stable, for any coherent quotient
sheaf $ \mathscr{Q} $ with $ \mathrm{rk}(\mathscr{Q}) > 0 $, one has
$ \lambda_E \leq  \lambda_{\sQ}.$ Hence,  \beq
\lambda_E=\lambda_E^-\eeq Fix $ P \in \mathrm{Herm}^+(E) $.  We set
$ h_0 = P $ on $ E $. For any $ \eps > 0 $, by Theorem \ref{main},
there exists a unique Hermitian metric $h_\eps$ on $E$ such that
\beq \label{epssolution} \Lambda_{\omega_g}\left(\sq
R^{h_\eps}\right)= (\lambda_E^- - \eps)h_\eps + \eps P. \eeq Taking
trace with respect to $ h_\eps $ in \eqref{epssolution} yields \beq
\label{approximateeps} \Lambda_{\omega_g}\left(\smo
\Theta^{h_\eps}\right) = (\lambda_E^--\eps)\mathrm{Id}_E + \eps
H_\eps^{-1}, \eeq where $ H_\eps = h_\eps \cdot h_0^{-1} \in
\Gamma(M,E^*\otimes E) $. On the other hand, we know \beq \int_M
\mathrm{tr}_E\left(\Lambda_{\omega_g}\left(\smo
\Theta^{h_\eps}\right) \right) \omega_g^n = \lambda_E
\rk(E)\int_M\omega_g^n=\lambda_E^- \rk(E)\int_M\omega_g^n. \eeq
Taking the trace of \eqref{approximateeps} and integrating over $M$,
we deduce that \beq \label{L1estimate}\int_M
\mathrm{tr}_E\left(H_\eps^{-1}\right) \omega_g^n = \rk(E)\int_M
\omega_g^n. \eeq On the other hand, by a straightforward computation
\be
\Delta_{\mathbb{C}}\mathrm{tr}_E\left(H_\eps^{-1}\right) & = & -\mathrm{tr}_E\left(\smo\Lambda_{\omega_g}\bar\p \left(\p^{h_0}H_\eps^{-1}\right)\right) = \mathrm{tr}_E\left(\smo\Lambda_{\omega_g}\bar\p \left(H_\eps^{-1} \cdot \p^{h_0}H_\eps \cdot H_\eps^{-1}\right)\right) \\
& = & -\mathrm{tr}_E\left(\smo\Lambda_{\omega_g}\left(H_\eps^{-1} \cdot \bar\p H_\eps \cdot H_\eps^{-1} \cdot \p^{h_0}H_\eps \cdot H_\eps^{-1}\right)\right) \\
&& + \mathrm{tr}_E\left(H_\eps^{-1} \cdot
\smo\Lambda_{\omega_g}\bar\p \left(\p^{h_0}H_\eps \cdot
H_\eps^{-1}\right)\right). \ee Since $\bp H_\eps$ and $\p^{h_0}
H_\eps$ are $h_0$-adjoint, by  (3) of Lemma \ref{linearalgebra}, we
have \beq
-\mathrm{tr}_E\left(\smo\Lambda_{\omega_g}\left(H_\eps^{-1} \cdot
\bar\p H_\eps \cdot H_\eps^{-1} \cdot \p^{h_0}H_\eps \cdot
H_\eps^{-1}\right)\right) \geq 0. \eeq By \eqref{conformalchange1}
and \eqref{approximateeps}, \be \mathrm{tr}_E\left(H_\eps^{-1} \cdot
\smo\Lambda_{\omega_g}\bar\p \left(\p^{h_0}H_\eps \cdot
H_\eps^{-1}\right)\right) & = & \mathrm{tr}_E\left(H_\eps^{-1} \cdot
\left(\Lambda_{\omega_g}\left(\smo \Theta^{h_\eps}\right) -
\Lambda_{\omega_g}\left(\smo
\Theta^{h_0}\right)\right)\right) \\
& = & \mathrm{tr}_E\left(H_\eps^{-1} \cdot
\left((\lambda_E^--\eps)\mathrm{Id}_E + \eps H_\eps^{-1} -
\Lambda_{\omega_g}\left(\smo
\Theta^{h_0}\right)\right)\right) \\
& \geq & (C-\eps) \cdot \mathrm{tr}_E(H_\eps^{-1}), \ee where $ C =
C(M,\omega_g,h_0,\lambda_E^-)$. Hence, for any $ 0 < \eps \leq 1 $,
one has \beq
\label{Laplaceestimate}\Delta_{\mathbb{C}}\mathrm{tr}_E\left(H_\eps^{-1}\right)
\geq (C-1) \cdot \mathrm{tr}_E\left(H_\eps^{-1}\right). \eeq
The Moser iteration shows
\beq \|\mathrm{tr}_E\bigl(H_\varepsilon^{-1}\bigr)\|^2_{L^\infty} \leq \tilde C_1 \|\mathrm{tr}_E\left(H_\varepsilon^{-1}\right)\|^2_{L^2}\eeq
where \(\tilde C_1 = \tilde C_1(M, \omega_g, C) > 0\). By \eqref{L1estimate},
one concludes that for all \(0 < \varepsilon \le 1\), \beq
\mathrm{tr}_E\bigl(H_\varepsilon^{-1}\bigr) \le C_1 \eeq where
$C_1=\tilde C_1  \rk(E)\int_M \omega_g^n$.
In particular,  one obtains the uniform $ C^0 $-estimate:
\beq \label{LowerC0estimateSS} H_\eps^{-1} \leq C_1\mathrm{Id}_E.
\eeq Moreover, for any $ 0 < \eps \leq 1 $, one has \beq
-C_2\mathrm{Id}_E \leq H_\eps^{-1} - \mathrm{Id}_E \leq
C_2\mathrm{Id}_E. \eeq where $ C_2 = \max\{C_1,1\}$.  In particular,
by \eqref{approximateeps}, \beq \lim_{\eps \> 0^+}
\nm{\Lambda_{\omega_g}\left(\sq \Theta^{h_\eps}\right) - \lambda_E^-
	\mathrm{Id}_E}_{L^\infty(M,\omega_g, h_0)} = \lim_{\eps \> 0^+}
\eps\nm{H_\eps^{-1} - \mathrm{Id}_E}_{L^\infty(M,\omega_g, h_0)} =
0. \eeq
Therefore, $ E $ is approximately Hermitian-Einstein.\\

Furthermore, we assume that $E$ is stable.  Under this assumption,  we established that,  for any $ 0 < \eps \leq 1 $,
\beq H_\eps^{-1} \leq C_1\mathrm{Id}_E. \eeq Now we turn to the
uniform upper bound for $h_\eps$. Let  \beq \Lambda_\eps := \sup_M
\lambda_{\max}(H_\eps)\eeq denote the supremum over $M$ of the
largest eigenvalue of $H_\eps$.  Define \beq \tilde h_\eps:=
\Lambda_\eps^{-1}\cdot h_\eps, \quad \tilde
H_\eps:=\Lambda_\eps^{-1}\cdot H_\eps. \eeq
Then it is clear that $ \tilde H_\eps\leq \mathrm{Id}_E$ with respect  to both $h_0$ and $h_\eps$.\\

Suppose that \beq \limsup_{\eps\>0+}\Lambda_\eps = \limsup_{\eps\>0+}\sup_M \lambda_{\max}(H_\eps)=+\infty.\eeq
By \eqref{L1estimate}, one has \beq \rk(E)\Lambda_\eps\int_M
\omega_g^n = \int_M \mathrm{tr}_E\left(\Lambda_\eps \cdot
H_\eps^{-1}\right) \omega_g^n = \int_M \mathrm{tr}_E\left(\tilde
H_\eps^{-1}\right) \omega_g^n \leq  \rk(E)\sup_M
\lambda_{\mathrm{max}}(\tilde H_\eps^{-1})\int_M \omega_g^n. \eeq
Hence, we deduce that $ \Lambda_\eps \leq
\sup_M\lambda_{\mathrm{max}}(\tilde H_\eps^{-1})$ and so \beq
\label{minimaleigenvlueHeps} \limsup_{\eps\>0+}
\sup_M\lambda_{\mathrm{max}}(\tilde H_\eps^{-1}) = +\infty. \eeq By
\eqref{epssolution}, one has \beq \label{lambdasoltildeheps}
\Lambda_{\omega_g}\left(\sq \Theta^{\tilde h_\eps}\right)
=\Lambda_{\omega_g}\left(\sq \Theta^{h_\eps}\right) =
(\lambda_E^--\eps)\mathrm{Id}_E + \eps H_\eps^{-1}.\eeq Hence, there
exists a constant $ C_3 = C_3(C_1,\lambda_E^-) > 0 $ such that for
any $ 0 < \eps \leq 1 $, \beq -C_3\mathrm{Id}_E \leq
\Lambda_{\omega_g}\left(\sq \Theta^{\tilde h_\eps}\right)  \leq
C_3\mathrm{Id}_E, \quad 0 < \tilde H_\eps \leq \mathrm{Id}_E. \eeq
By using a strategy similar to that in the proof of
Theorem~\ref{main}, there exist a subsequence $\{\tilde
H_{\eps_i}\}$ of $\tilde H_\eps$, and a sequence $\{\tau_i\}\subset
(0,1/2]$ decreasing to $0$ such that

\bd\item $\lim_i\Lambda_{\eps_i}=+\infty$.
\item $ \{\tilde H_{\eps_i}\}$ converges to  some  $\tilde H_0$ in weak $ W^{1,2}(M, \omega_g, h_0) $ sense.
\item$ \{\tilde H^{2\tau_i}_0\}$ converges to some $ \tilde H $ in weak $ W^{1,2}(M, \omega_g, h_0) $ sense.
\item$ \pi_1 = \mathrm{Id}_E - \tilde H $ defines a weak holomorphic projection with respect to $ h_0 $. Let $\cF_1$ be the coherent subsheaf of $E$ defined by $\pi_1$. By \eqref{minimaleigenvlueHeps} and a uniform positive lower bound of $\left\|\tilde H_\eps\right\|_{L^2;h_0}$,  one has
\beq 0 < \mathrm{rk}(\mathscr{F}_1) < \mathrm{rk}(E).\eeq \ed
Moreover, denoting $\mathscr{Q}_1 = E/\mathcal{F}_1$,  we obtain the
following estimate by a proof similar to that of
\eqref{quotientsheafQ}: \beq \label{quotientsheafQ1}
\mathrm{deg}_{\omega_g}(\mathscr{Q}_1) \leq \limsup_{\tau_j \>
	0^+}\limsup_{\eps_i \> 0^+}\int_M\tr_E\left(\Lambda_{\om_g}\left(\sq
\Theta^{h_{\eps_i}}\right)\cdot \tilde
H_{\eps_i}^{2\tau_j}\right)\om_g^n. \eeq On the other hand, by
\eqref{lambdasoltildeheps}, for any $ 2\tau_j < 1 $,
\be && \limsup_{\eps_i \> 0^+} \left(\Lambda_{\om_g}\left(\sq \Theta^{h_{\eps_i}}\right), \tilde H_{\eps_i}^{2\tau_j}\right)_{h_0} = \limsup_{\eps_i \> 0^+} \left((\lambda_E^--\eps_i)\mathrm{Id}_E + \eps_i H_{\eps_i}^{-1}, \tilde H_{\eps_i}^{2\tau_j}\right)_{h_0} \\
& \leq & \limsup_{\eps_i \> 0^+} \left((\lambda_E^- -\eps_i)\mathrm{Id}_E, \tilde H_{\eps_i}^{2\tau_k}\right)_{h_0} + \limsup_{\eps_i \> 0^+} \eps_i\left(H_{\eps_i}^{-1}, \tilde H_{\eps_i}^{2\tau_j}\right)_{h_0} \\
& = & \lambda_E^-\left(\mathrm{Id}_E,\tilde
H^{2\tau_j}_0\right)_{h_0}. \ee As $ \tau_j \> 0^+ $, we deduce that
\be \limsup_{\tau_j \> 0^+}\limsup_{\eps_i \> 0^+}
\left(\Lambda_{\om_g}\left(\sq \Theta^{h_{\eps_i}}\right), \tilde
H_{\eps_i}^{2\tau_j}\right)_{h_0}
& \leq & \limsup_{\tau_j \> 0^+}\lambda_E^- \left(\mathrm{Id}_E,\tilde H^{2\tau_j}_0\right)_{h_0} \\
& = & \lambda_E^-\left(\mathrm{Id}_E, \tilde H \right)_{h_0} =
\lambda_E^-\mathrm{rk}(\mathscr{Q}_1)\left(\int_M \omega_g^n\right).
\ee Hence, by \eqref{quotientsheafQ1}, one has \beq
\label{slopeQ1estimate}  \mathrm{deg}_{\omega_g}(\mathscr{Q}_1) \leq
\lambda_E^-\mathrm{rk}(\mathscr{Q}_1)\left(\int_M \omega_g^n\right).
\eeq By \eqref{slopeQ1estimate} and the fact that $
\lambda_E^-=\lambda_E $, one has \beq
\frac{\mathrm{deg}_{\omega_g}(E)}{\mathrm{rk}(E)} \geq
\frac{\mathrm{deg}_{\omega_g}(\mathscr{Q}_1)}{\mathrm{rk}(\mathscr{Q}_1)}.
\eeq
Moreover, $ 0 < \mathrm{rk}(\mathscr{F}_1) < \mathrm{rk}(E) $,  and $ 0 < \mathrm{rk}(\mathscr{Q}_1) < \mathrm{rk}(E) $. This contradicts the assumption that $ E $ is stable. \\

Hence, there exists a uniform constant $ C > 0 $ such that \beq
\limsup_{\eps\>0+}\Lambda_\eps = \limsup_{\eps\>0+}\sup_M
\lambda_{\max}(H_\eps)=C<+\infty.\eeq Then we have \beq h_\eps \leq
C h_0, \eeq for all  $ \eps$ close to $0$. We rewrite the equation
\eqref{epssolution} as \beq  \Lambda_{\omega_g}\left(\sq
R^{h_\eps}\right)= \lambda_\eps h_\eps + P_\eps, \eeq where
$\lambda_\eps=\lambda_E^- - \eps$ and $P_\eps=\eps P$.  It is
obvious that $\{P_\eps\}$ converges smoothly to $P=0$. By Theorem
\ref{C^0toC^infty}, there exists a subsequence $\{h_{\eps_i}\}$
which converges smoothly to a smooth Hermitian metric $h$ on $E$,
and this limit satisfies
\begin{equation}
\Lambda_{\omega_g}\left(\sq  R^{h} \right)= \lambda_E^- h .
\end{equation}
In particular, $(E,h)$ is Hermitian-Einstein. This completes the
proof. \eproof \bremark If $E$ is not stable, Theorem~\ref{main}
guarantees the existence of a family of Hermitian metrics
$h_\varepsilon$ solving \beq \Lambda_{\omega_g}\left(\sqrt{-1}
R^{h_\varepsilon}\right) = (\lambda_E^- -
\varepsilon)\,h_\varepsilon + \varepsilon P, \eeq for any fixed
$P\in\mathrm{Herm}^+(E)$  and $\varepsilon>0$.  It is natural to study the
limit of $h_\eps$ as $\varepsilon\to 0^+$. \eremark

\vskip 1\baselineskip

\end{document}